\documentclass[11pt]{article}

\usepackage{amssymb}
\usepackage{amsmath}
\usepackage{amsthm}
\usepackage{graphicx}

\newcommand{\R}{\mathbb R}
\newcommand{\C}{\mathbb C}
\newcommand{\V}{\mathbb V}
\newcommand{\St}{\mathbb S}
\newcommand{\Oh}{{\cal O}}
\newcommand{\st}{{\triangle t}}
\newcommand{\hQ}{{\hat{Q}}}
\newcommand{\hR}{{\hat{R}}}
\renewcommand{\i}{\mathrm i}
\renewcommand{\Im}{\mathrm{Im}}
\renewcommand{\Re}{\mathrm{Re}}

\newtheorem{remark}{Remark}[section]
\newtheorem{theorem}{Theorem}[section]
\newtheorem{lemma}{Lemma}[section]
\newtheorem{example}{Example}[section]

\author{J.~Lang, W.~Hundsdorfer}
\title{Extrapolation-Based Implicit-Explicit Peer Methods with Optimised
Stability Regions}
\author{
Jens Lang\footnote{corresponding author}\\
{\small \it Technische Universit\"at Darmstadt} \\ {\small \it
Dolivostra{\ss}e 15, 64293 Darmstadt, Germany}\\
{\small lang@mathematik.tu-darmstadt.de} \\ \\
Willem Hundsdorfer \\
{\small \it Center for Mathematics and Computer Science} \\
{\small \it P.O.\ Box 94079, 1090 GB Amsterdam, The Netherlands}\\
{\small Willem.Hundsdorfer@cwi.nl}}
\date{January 19, 2017}
\begin{document}
\maketitle

\begin{abstract}
In this paper we investigate a new class of implicit-explicit (IMEX)
two-step methods of Peer type for systems of ordinary differential equations
with both non-stiff and stiff parts included in the source term. An
extrapolation approach based on already computed stage values
is applied to construct IMEX methods with favourable
stability properties. Optimised IMEX-Peer methods of
order $p=2,3,4$, are given as result of a search algorithm carefully designed
to balance the size of the stability regions and the extrapolation errors.
Numerical experiments and a comparison to other implicit-explicit methods
are included.
\end{abstract}


\section{Introduction}
Many initial value problems arising in practice are in a form
$u' = F_0(u) + F_1(u)$, where $F_0$ is a non-stiff or mildly stiff
part and $F_1$ is a stiff contribution. Implicit-explicit (IMEX)
methods use this decomposition by treating only the $F_1$
contribution in an implicit fashion. The advantage of lower costs
for explicit schemes is combined with the favourable stability
properties of implicit schemes to enhance the overall computational
efficiency.

In this paper we will consider IMEX methods based on implicit
Peer methods. These methods introduced by Schmitt, Weiner
and co-workers
\cite{BeckWeinerPodhaiskySchmitt2012,SchmittWeiner2004,SchmittWeinerErdmann2005}
are a very comprehensive class of general linear methods
(GLMs) in which the approximations in all stages have the same order.
Peer methods can be viewed as a natural generalisation of linear
multistep methods in the sense that each of the stages is a linear
multistep method itself. Due to their multi-stage structure they inherit
good stability properties and an easy step size change in every time step
from one-step methods without suffering from order reduction for stiff problems.
The property that the peer stage values have the same order of accuracy can
be conveniently exploited to construct related explicit methods by using
extrapolation. The combination of these implicit and explicit methods leads
in a natural way to IMEX methods with the same order as the original
implicit method. This idea was first used by Crouzeix \cite{Crouzeix1980}
with linear multistep methods of BDF type.
IMEX-Peer methods are competitive alternatives to classic IMEX methods
for large stiff problems. Higher-order IMEX Runge-Kutta methods are known to
suffer from possible order reduction and serious efficiency loss for stiff
problems. Moreover, the increasing number of necessary coupling conditions
makes their construction difficult.

Recently, the same extrapolation approach was used by
Cardone, Jackiewicz, Sandu and Zhang \cite{CardoneJackiewiczSanduZhang2014}
starting with diagonally implicit multistage integration methods (DIMSIMs).
In those general linear methods
the implicit internal stages are followed by explicit stages. Due to
these explicit stages the linear stability properties of the resulting
IMEX methods are less favourable than for the IMEX-Peer methods. Higher-order
IMEX-GLMs were constructed by Zhang, Sandu and Blaise
\cite{ZhangSanduBlaise2016}, based on an earlier developed partitioned GLM
framework of the same authors in \cite{ZhangSanduBlaise2014}.
Bra\'s, Izzo and Jackiewicz investigated IMEX-GLMs of order
up to four with inherent Runge-Kutta stability \cite{BrasIzzoJackiewicz2016}.

In Section 2 of this paper we present the framework to obtain
IMEX-Peer methods based on suitable implicit methods. The construction
of specific classes of methods is performed in Section 3. Along with
IMEX-BDF methods, which also fit in the Peer form, we will construct
IMEX-Peer methods based on the implicit methods of order 3 and 4 that
were developed by Beck, Weiner, Podhaisky and Schmitt
\cite{BeckWeinerPodhaiskySchmitt2012}. Comparison
of the stability regions of the methods shows promising results for
the latter methods. This is confirmed in the numerical experiments in
Section 5 for two advection-reaction problems with stiff reactions
and a reaction-diffusion problem, where the diffusion leads to stiffness.

\section{Implicit-Explicit Peer Methods Based on Extrapolation}
\subsection{Implicit Peer methods}
To solve initial value problems in the vector space $\V=\R^m, m\ge 1$,
\begin{equation}
u'(t) = F(u(t)),\quad u(0)=u_0\in\V\,,
\end{equation}
we consider the so-called Peer methods introduced by Schmitt, Weiner and co-workers
\cite{SchmittWeiner2004,BeckWeinerPodhaiskySchmitt2012}.
An $s$-stage Peer method provides approximations
\begin{equation}
w_n = [w_{n,1},\ldots,w_{n,s}]^T\in\V^s,\quad w_{n,i}\approx u(t_n+c_i\st)\,,
\end{equation}
where $t_n = n \st$, $n\ge0$, and the nodes $c_i\in\R$ are
such that $c_i\ne c_j$ if $i\ne j$, and $c_s=1$. The starting
vector $w_0=[w_{0,i}]\in\V^s$ is supposed to be given, or computed from a Runge-Kutta
method, for example.

Peer methods are general linear methods, based on the requirement that all
approximations $w_{n,j}$ have the same order. Here, we will primarily aim at order
$p=s$. With $s\times s$ coefficient matrices $P=(p_{ij})$, $Q=(q_{ij})$,
$R=(r_{ij})$, and the $m\times m$ identity matrix $I$, the usual general form of
the implicit methods of this Peer type is
\begin{equation}
\label{implpeer1}
w_n = (P\otimes I)w_{n-1} + \st (Q\otimes I)F(w_{n-1}) + \st (R\otimes I)F(w_n)\,.
\end{equation}
where $F(w)=[F(w_i)]\in\V^s$ is the application of $F$ to all components
of $w\in\V^s$. In the following, for an $s\times s$-matrix we will use the same
symbol for its Kronecker product
with the identity matrix as mapping from the space $\V^s$ to itself.
Then, (\ref{implpeer1}) simply reads
\begin{equation}
\label{implpeer}
w_n = Pw_{n-1} + \st QF(w_{n-1}) + \st RF(w_n)\,.
\end{equation}
The matrix $R$ is taken to be lower triangular, giving diagonally implicit methods,
with diagonal $R$ if parallelism is a special case of interest
\cite{SchmittWeinerErdmann2005}. Implicit peer methods with good stability properties,
i.e., $L(\alpha)$-stability with large angles $\alpha$, can be found by taking
$Q=0$ \cite{BeckWeinerPodhaiskySchmitt2012}. We will choose these methods to construct
implicit-explicit peer methods based on extrapolation. Then the method reads
\begin{equation}
\label{implpeerq0}
w_n = P w_{n-1} + \st RF(w_n)\,.
\end{equation}
Some requirements or desirable properties are briefly discussed here for the implicit
method (\ref{implpeerq0}).\\

\noindent {\bf Zero-stability.} The matrix $P$ should be power bounded to have stability
for the trivial problem $u'(t)=0$. Let $\text{spr}(P)$ be the spectral radius of $P$.
Since one eigenvalue of $P$ will be equal to $1$ for pre-consistency, the requirement
of zero-stability means
\begin{equation}
\text{spr}(P)=1 \; \text{and eigenvalues with modulus 1 are not defective}.
\end{equation}
This requirement was enforced by Schmitt, Weiner et al. by taking $P$ such that
one eigenvalue equals $1$ and the others are $0$. This choice, called optimal zero-stability,
made the construction of methods more tractable. We will also look at methods that are
strongly zero-stable, where $P$ has one eigenvalue $1$ and the other eigenvalues have
modulus less than $1$. This holds for example for the well-known BDF methods.\\

\noindent {\bf Accuracy.} Let $e=(1,\ldots,1)^T\in\R^s$. It will be assumed that
\begin{equation}
\label{precons}
P e = e\,.
\end{equation}
This is the so-called
pre-consistency condition, which means that for the trivial equation $u'(t)=0$,
we get solutions $w_{n,i}=1$ provided that $w_{0,j}=1,\;j=1,\ldots,s$. Inserting exact solution
values $w(t_n)=[u(t_n+c_i\st)]\in\V^s$ in the implicit scheme (\ref{implpeerq0})
gives the residual-type local errors
\begin{equation}
r_n = w(t_n) - P w(t_{n-1}) - \st R w'(t_n)\,.
\end{equation}
Let $c=(c_1,\ldots,c_s)^T$ with
point-wise powers $c^j=(c_1^j,\ldots,c_s^j)^T$. Then Taylor expansion gives
\begin{eqnarray}
w(t_n) &\!=\!& e\otimes u(t_n) + \st c\otimes u'(t_n) +
  \frac{1}{2}\st^2 c^2\otimes u''(t_n) + \ldots\,\\
w(t_{n-1}) &\!=\!& e\otimes u(t_n) + \st (c-e)\otimes u'(t_n) +
  \frac{1}{2}\st^2 (c-e)^2\otimes u''(t_n) + \ldots,
\end{eqnarray}
from which we obtain
\begin{equation}
r_n = \sum_{j\ge 1} \st^j d_j\otimes u^{(j)}(t_n)
\end{equation}
with
\begin{equation}
d_j = \frac{1}{j!}\left( c^j - P(c-e)^j - jRc^{j-1}\right)\,.
\end{equation}
The method is said to have (stage) order $q$ if (\ref{precons}) holds and
$d_j=0$ for $j=1,2,\ldots,q$. We will be interested in methods with
(stage) order $s$. With the Vandermonde matrices
\begin{equation}
\label{eq:V0V1}
V_0 = \big(c_i^{j-1}\big), \qquad
V_1 = \big((c_i-1)^{j-1}\big),\qquad i,j=1,\ldots,s,
\end{equation}
and the diagonal matrices $C=\text{diag}(c_1,c_2,\ldots,c_s)$,
$ D=\text{diag}(1,2,\ldots,s)$, the conditions for having stage order
$s$ with the implicit method (\ref{implpeerq0}) are
\begin{equation}
CV_0 - P(C - I)V_1 - RV_0 D = 0\,.
\end{equation}

\begin{remark}[superconvergence]
For a method with stage order $q$, it is possible to have convergence with
order equal to $q+1$. This is discussed under the heading super-convergence
in the book of Strehmel, Weiner and Podhaisky
\cite[Sect.\,5.3]{StrehmelWeinerPodhaisky2012}
for non-stiff problems.
It is related to the definition of order of consistency for general
linear methods as given in
\cite[Sect.\,III.8]{HairerNoersettWanner1993}.
Similar results for stiff systems can be found in \cite{Hundsdorfer1994}.
\end{remark}

\subsection{Extrapolation}
Based on an implicit method with order $s$, a related explicit method can
be found by extrapolation, leading to implicit-explicit methods.
This is a well-known procedure for linear multistep methods, see for
instance Crouzeix \cite{Crouzeix1980} or the review in the book of
Hundsdorfer and Verwer \cite[Sect.\,IV.4.2]{HundsdorferVerwer2003}.
Recently this idea was also used with a class of general linear methods,
the so-called diagonally implicit multistage integration methods (DIMSIMs),
by Cardone, Jackiewicz, Sandu, and Zhang \cite{CardoneJackiewiczSanduZhang2014}.
Here, we will use this extrapolation idea to obtain implicit-explicit Peer
methods.

Having an implicit method, where all approximations $w_{n,j}$ have order $s$,
we can obtain a corresponding explicit method by extrapolation using a
Lagrange polynomial of degree $s-1$, giving
$\varphi(t_{n,i})=\sum_j s_{ij}\varphi(t_{n-1,j})+\Oh(\st^s)$
for smooth functions $\varphi$, with $t_{n,i}=t_n+c_i\st$.
The extrapolation coefficients are given by
$s_{ij}=\Pi_{k\ne j}(c_i-c_k+1)/(c_j-c_k)$.

We can apply this extrapolation with $\varphi(t)=F(u(t))$. Starting from
the implicit method (\ref{implpeerq0}), this yields the explicit method
\begin{equation}
\label{explpeer}
w_n = P w_{n-1} + \st\hQ F(w_{n-1})\,,
\end{equation}
with coefficient matrix $\hQ=(\hat{q}_{ij})$ given by $\hQ=RS$, where $S=(s_{ij})$.
By the construction, all the stages have again order $s$, at least, so (\ref{explpeer})
is an explicit Peer method.

The extrapolation may be improved by using the last available information, whereby
a value $\varphi(t_{n,i})$ is found as linear combination of some of the values
$\varphi(t_{n-1,j})$ together with the most recent values $\varphi(t_{n,j})$,
$j=1,\ldots,i-1$, say
\begin{equation}
\varphi(t_{n,i})=\sum_{j} s^{(1)}_{ij}\varphi(t_{n-1,j}) +
\sum_{j\le i-1}s^{(2)}_{ij}\varphi(t_{n,j})
+ \Oh(\st^s)\,,\quad i=1,\ldots,s.
\end{equation}
Setting $S_1=(s^{(1)}_{ij})$, $S_2=(s^{(2)}_{ij})$, this will lead to an explicit Peer
method of the form
\begin{equation}
\label{explpeerexpol}
w_n = P w_{n-1} + \st\hQ F(w_{n-1}) + \st\hR F(w_n)
\end{equation}
with
\begin{equation}
\hQ = RS_1,\quad \hR = RS_2\,.
\end{equation}
Note that $\hR$ is strictly lower triangular, since $R$ is lower triangular and $S_2$
is strictly lower triangular.

Defining vectors $\Phi_m=[\varphi(t_{m,i})]\in\V^s$, the error vector for the
extrapolation, $\delta_n=\Phi_n-S_1\Phi_{n-1}- S_2\Phi_n$, can
be expanded in a Taylor series at $t_n$,
\begin{equation}
\label{errexpol}
\delta_n = (I-S_1-S_2)e\otimes\varphi(t_n) + \sum_{j\ge 1}\frac{1}{j!}
\left( (I-S_2)c^j - S_1(c-e)^j\right)\otimes\varphi^{(j)}(t_n)\st^j\,.
\end{equation}
Therefore, the conditions for stage order $s$ read
\begin{equation}
\label{eq:stage_order_s}
(I-S_2) c^j - S_1 (c-e)^j = 0,\quad 0\le j\le s-1\,,
\end{equation}
which is equivalent to the relation $S_1V_1=(I-S_2)V_0$. The choice of a strictly
lower triangular $S_2$ thus determines $S_1$.

\subsection{Implicit-Explicit Peer Methods}
The combination of the related implicit and explicit methods (\ref{implpeer}),
(\ref{explpeerexpol}) can now be used to construct an implicit-explicit (IMEX)
method for systems of the form
\begin{equation}
\label{ode}
u'(t) = F_0(u(t)) + F_1(u(t))\,,
\end{equation}
where $F_0$ will represent the non-stiff or mildly stiff part, and $F_1$ gives
the stiff part of the equation. The resulting IMEX scheme is
\begin{equation}
\label{imex-peer}
w_n = P w_{n-1} + \st\hQ F_0(w_{n-1}) + \st\hR F_0(w_{n})
+ \st R F_1(w_{n})\,.
\end{equation}
The extrapolation idea is used here only on the $F_0$. For non-stiff problems,
this IMEX method will have order $s$ for any decomposition $F=F_0+F_1$. However,
for stiff problems it should be required that the derivatives of
$\varphi_k(t)=F_k(u(t)), \;k=0,1,$ are bounded by a moderate constant which is
not affected by the stiffness parameters, such as the spatial mesh width $h$ for
semi-discrete systems obtained from PDEs.

\begin{remark}[linearly implicit methods]
If $J\approx F'(u)$, such that $F(u)-Ju$ is a non-stiff or mildly stiff term,
we can consider the decomposition $F_0 (u)=F(u)-Ju$ and $F_1(u)=Ju$ as
special case of (\ref{ode}). This gives the linearly implicit Peer method
\begin{equation}
w_n = P w_{n-1} + \st\hQ F(w_{n-1}) + \st\hR F(w_{n})
- \st\hQ J w_{n-1}  + \st(R-\hR) J w_{n}\,.
\end{equation}
By the above construction, leading to the IMEX scheme (\ref{imex-peer}), all stages
will be consistent of order $s$. However, since $F_1$ is linear here, it is possible
that the order conditions (\ref{eq:stage_order_s}) are in fact a bit too strong.
\end{remark}

The standard local consistency analysis for the IMEX-Peer method (\ref{imex-peer})
with exact solution values $u(t_{n,i})$ yields for the residual-type local errors
\begin{equation}
r_n = E_{im} + \st R\,E_{ex} + \Oh\left(\st^{s+2}\right)\,,
\end{equation}
where $E_{im} = \st^{s+1} d_{s+1}\otimes u^{(s+1)}(t_n)$ is the leading
error term of the corresponding implicit Peer method with constant time steps.
Replacing in (\ref{errexpol}) $\varphi(t)$ by $F_0(u(t))$ taken as function
of $t$, we find for the leading error vector of the extrapolation,
\begin{equation}
E_{ex} = \frac{\st^s}{s!}\big( (I-S_2)c^s - S_1(c-e)^s\big)\otimes
\frac{d^s}{dt^s}F_0(u(t_n))\,.
\end{equation}
Together with zero-stability of the implicit Peer method and standard
convergence arguments, we have the following result for the IMEX scheme
(\ref{imex-peer}) applied to non-stiff problems:
\begin{theorem}
Let the $s$-stage implicit Peer method (\ref{implpeerq0}) with
coefficients $(c,P,R)$ be zero-stable and suppose its stage order is
equal to $s$.
Let the starting values satisfy $w_{0,i}-u(t_0+c_i\st)=\Oh(\st^s),\,
i=1,\ldots,s$.
Then the IMEX scheme (\ref{imex-peer}) with $\hR=RS_2$ and
$\hQ=R(I-S_2)V_0V_1^{-1}$ is convergent of order $s$ for constant step size
and arbitrary strictly lower triangular matrix $S_2$.
\end{theorem}
Note that $E_{im}$ and $E_{ex}$ are not influenced by stiffness, and the same
result will therefore hold for stiff problems provided suitable linear or
nonlinear stability conditions are satisfied.

For later use we define the following two error constants:
\begin{equation}
c_{im} = \big\|d_{s+1}\big\| = \frac{1}{(s+1)!} \, \big\|
\left( c^{s+1}-P(c-e)^{s+1}-(s+1)Rc^s\right) \big\|
\end{equation}
and
\begin{equation}
\label{c_ex}
c_{ex} = \frac{1}{s!} \, \big\| \big( (R-\hR)c^s - \hQ(c-e)^s\big) \big\|\,.
\end{equation}
with $\|\cdot\|$ being the Euclidean norm in $\R^s$.
The first one is the error constant of the implicit Peer method and the second one is
related to the extrapolation process.

\subsection{Stability of IMEX-Peer Methods}
We consider the general test equation
\begin{equation}
\label{imex-test-eq}
y'(t) = \lambda_0 y(t) + \lambda_1y(t),\quad t\ge 0,
\end{equation}
with complex parameters $\lambda_0$ and $\lambda_1$.
Define $z_i=h\lambda_i, i=0,1$.
Applying an IMEX-Peer method to (\ref{imex-test-eq}) gives
\begin{equation}
\label{imex-peer-lin}
w_{n+1} = (I-z_0RS_2-z_1R)^{-1}(P+z_0RS_1)w_n\,.
\end{equation}
This can be compactly written as $w_{n+1}=M(z_0,z_1)w_n$.
For given $z_0$ and $z_1$, stability is ensured if
\begin{equation}
\label{imex-peer-stabmat}
\rho (M(z_0,z_1)) < 1.
\end{equation}
The stability function of the IMEX-Peer method is defined as the
characteristic polynomial of the stability matrix $M(z_0,z_1)$:
\begin{equation}
\label{polM}
\zeta(w,z_0,z_1) = \det (wI - M(z_0,z_1))\,.
\end{equation}
Consequently, the IMEX-Peer method is stable for given
$z_0,z_1\in\C$ if all the roots $w_i(z_0,z_1), i=1,\ldots,s$,
of the stability function $\zeta(w,z_0,z_1)$ are inside the unit
circle.

The higher order implicit Peer methods considered here
are $L(\alpha)$-stable with respect to the implicit part $z_1\in\C$.
Therefore, we introduce for $\beta \in [0,\frac{1}{2}\pi]$ the sets
\begin{equation}
\St_\beta = \{z_0\in\C: (\ref{imex-peer-stabmat})\text{ holds for any }
z_1 \in\C \text{ with }|\Im (z_1)|\le -\tan(\beta)\cdot\Re(z_1) \}
\end{equation}
in the left-half complex plane.
In order to compute these sets for specific angles $\beta$ with
$0\le\beta\le\alpha$, we first define for fixed $y\in\R$,
\begin{equation}
\St_{\beta,y} =
\left\{
z_0\in\C: (\ref{imex-peer-stabmat})\text{ holds for fixed }
z_1=-|y|/\tan(\beta) + y\,\i
\right\}
\end{equation}
and find then $\St_\beta$ from the intersection of all
$\St_{\beta,y},\; y\in\R$, which follows from the maximum principle.
The set $\St_E:=\St_{\beta,0}$ is independent of $\beta$ and
corresponds to the stability region of the explicit method.
Since $\St_\alpha\subset \St_E$, the goal is to construct IMEX-Peer methods
for which $\St_E$ is large and $\St_E\backslash \St_\alpha$ is as small as possible for
angles $\alpha$ that are close to $\pi/2$, whereas the error
constant for the extrapolation, $c_{ex}$, is still of moderate size.

The boundary locus method can be used to compute the boundary of
$\St_{\beta,y}$:
\begin{equation}
\partial \St_{\beta,y} = \{ z_0\in\C:
\zeta(e^{\theta\,\i},z_0,-|y|/\tan(\beta)+y\,\i)=0,\;
\theta\in [0,2\,\pi) \}\,.
\end{equation}
Varying the eigenvalues $w=e^{\theta\,\i}$ for
fixed $z_1(y)=-|y|/\tan(\beta)+y\,\i, y\in\R,$ allows to reformulate
the eigenvalue problem $M(z_0,z_1)x=wx$ into an eigenvalue problem
for $z_0(e^{\theta\,\i},y)$, i.e., $G(e^{\theta\,\i},z_1(y))x=z_0x$ with
\begin{equation}
G(e^{\theta\,\i},z_1(y)) = (e^{\theta\,\i}\hat{R} +
\hat{Q})^{-1}(e^{\theta\,\i}I-e^{\theta\,\i}z_1(y)R-P)\,.
\end{equation}
The set of all eigenvalues $z_0(e^{\theta\,\i},y)$ contains the boundary
of $\St_{\beta,y}$.

In order to approximate the boundary of the stability region $\St_\beta$,
we will follow the approach that was
successfully applied in \cite{CardoneJackiewiczSanduZhang2014} for DIMSIMs.
There the intersection point $z_0(y)\in\C$ of the boundary $\partial \St_{\beta,y}$ with
a ray $y_0=mx_0$ is computed by the bisection method with the termination condition
\begin{equation}
\label{cond_bi}
\left| \max_{i=1,\ldots,s}\left|w_i(z_0(y),-|y|/\tan(\beta)+y\,\i\right|-1 \right| \le tol
\end{equation}
for an appropriate accuracy tolerance $tol$. We set $tol=1e\!-\!5$ and start
with a large enough interval on the real axis so that (\ref{cond_bi}) is not
satisfied in the first iteration step. Then, for a fixed value $m\in\R$, the intersection
point $z_0(y)=x_0(y)+y_0(y)\,\i$ of the corresponding ray and the boundary $\partial \St_\beta$ is
determined by minimising $|x_0(y)|$ as function of $y$, which is passed to the {\sc Matlab}-routine
{\it fminsearch}. Finally, varying the parameter value $m$ delivers a polygonal
approximation to $\partial \St_\beta$, which is used to compute the size of the
area of the stability region, $|\St_\beta|$.

\section{Construction of IMEX-Peer Methods}
\subsection{IMEX-Peer Methods with Equidistant Nodes}
We will first consider IMEX Peer methods with equidistant nodes. A good
candidate within this class are the IMEX-BDF methods introduced in
\cite{Crouzeix1980,Varah1980}.
To have a closer resemblance with the usual Peer form, we formulate these
BDF methods with $s$ small steps of length $\st/s$.

In the following, let $a_0, a_1,\ldots,a_s$ be the coefficients of the
$s$-step BDF method; cf.\ Table~\ref{tab-bdf-coeff}.
Starting with approximate solutions
$u_{n-1+i/s}\approx u(t_{n-1} + \st\,i/s)$ for $i=1,\ldots,s$,
we have for $s$ steps, $k=1,\ldots,s$, with step-size $\st/s$ of the
IMEX-BDF method:
\begin{equation}
\label{imexbdf}
\sum_{i=0}^{s} a_i\,u_{n+(k-i)/s} =
\frac{\st}{s} \sum_{i=1}^{s} b_i\,F_0(u_{n-1+(k+i-1)/s})
+ \frac{\st}{s} F_1(u_{n+k/s}),
\end{equation}
where $(b_1,\ldots,b_s) = e_s^T S$, $e_s^T=(0,\ldots,0,1)$ and
$S=(s_{ij})=V_0V_1^{-1}$ is defined by (\ref{eq:V0V1}) with the
normalised vector $c=(0,1,\ldots,s-1)^T$.
In order to obtain the standard form of a Peer method, we set
\begin{equation}
w_{n-1} =
\begin{pmatrix}
u_{n-1+1/s}\\
u_{n-1+2/s}\\
\vdots\\
u_{n}
\end{pmatrix}
\qquad\mbox{and}\qquad
w_{n} =
\begin{pmatrix}
u_{n+1/s}\\
u_{n+2/s}\\
\vdots\\
u_{n+1}
\end{pmatrix}\,.
\end{equation}
This yields
\begin{equation}
A_2w_n = -A_1 w_{n-1} +\frac{\st}{s}B_1F_0(w_{n-1}) +
\frac{\st}{s}B_2F_0(w_{n}) + \frac{\st}{s}F_1(w_{n})
\end{equation}
with the matrices
\begin{equation}
A_1 =
\begin{pmatrix}
a_s & a_{s-1} & \cdots & a_1\\
  0 & a_s     & \cdots & a_2\\
\vdots & \vdots & \ddots & \vdots\\
  0 & 0 & \cdots & a_s
\end{pmatrix},\quad
A_2 =
\begin{pmatrix}
a_0 & 0 & \cdots & 0\\
a_1 & a_0 & \cdots & 0\\
\vdots & \vdots & \ddots & \vdots\\
a_{s-1} & a_{s-2} & \cdots & a_0
\end{pmatrix}
\end{equation}
and
\begin{equation}
B_1 =
\begin{pmatrix}
s_{s1} & s_{s2} & \cdots & s_{ss}\\
  0 & s_{s1}     & \cdots & s_{s,s-1}\\
\vdots & \vdots & \ddots & \vdots\\
  0 & 0 & \cdots & s_{s1}
\end{pmatrix},\quad
B_2 =
\begin{pmatrix}
0 & 0 & \cdots & 0\\
s_{ss} & 0 & \cdots & 0\\
\vdots & \ddots & \ddots & \vdots\\
s_{s2} & \cdots & s_{ss} & 0\\
\end{pmatrix}.\\[2mm]
\end{equation}
Note that $A_2$ is invertible since always $a_0\ne 0$. The coefficients
of the equivalent IMEX-Peer method are then given by the following Lemma:
\begin{lemma}
An $s$-stage IMEX-BDF method (\ref{imexbdf}) with $s$ steps of
length $\st/s$ is equivalent to an $s$-stage IMEX-Peer method (\ref{imex-peer})
with node vector $c=(1/s,2/s,\ldots,1)^T$ and coefficient matrices
\begin{equation}
P=-A_2^{-1}A_1,\;\hQ=(1/s)A_2^{-1}B_1,\;\hR=(1/s)A_2^{-1}B_2,
\text{ and } R=(1/s)A_2^{-1}\,.
\end{equation}
\end{lemma}
IMEX-BDF methods have proven to work very well and therefore they are a good
target for general IMEX-Peer methods with $p=s$.

\begin{example} Exemplarily, the (peer-)coefficients of the IMEX-BDF3 method
with three steps of length $\st/3$ are given. The node vector $c=(0,1,2)^T$
yields $e_3^TS=(1,-3,3)$. Thus, the matrices in (\ref{imex-peer}) are:
\begin{equation}
P =
\begin{pmatrix}
\frac{2}{11} & -\frac{9}{11} & \frac{18}{11}\\[2mm]
\frac{36}{121} & -\frac{140}{121} & \frac{225}{121}\\[2mm]
\frac{450}{1331} & -\frac{1629}{1331} & \frac{2510}{1331}
\end{pmatrix},\quad
R =
\begin{pmatrix}
\frac{2}{11} & 0 & 0\\[2mm]
\frac{36}{121} & \frac{2}{11} & 0\\[2mm]
\frac{450}{1331} & \frac{36}{121} & \frac{2}{11}
\end{pmatrix}\,,
\end{equation}
\begin{equation}
\hat{Q} =
\begin{pmatrix}
\frac{2}{11} & -\frac{6}{11} & \frac{6}{11}\\[2mm]
\frac{36}{121} & -\frac{86}{121} & \frac{42}{121}\\[2mm]
\frac{450}{1331} & -\frac{954}{1331} & \frac{404}{1331}
\end{pmatrix},\quad
\hat{R} =
\begin{pmatrix}
0 & 0 & 0\\[2mm]
\frac{6}{11} & 0 & 0\\[2mm]
\frac{42}{121} & \frac{6}{11} & 0
\end{pmatrix}\,.\\[2mm]
\end{equation}
The eigenvalues of $P$ are $1$ and $(-119\pm 27\sqrt{39}i)/2662$.
\end{example}

\begin{table}[t]
\centering
\begin{tabular}{|r|rrrrr|c|}
\hline
$s$ & $a_0$ & $a_1$ & $a_2$ & $a_3$ & $a_4$ & $\alpha$\\
\hline\rule{0mm}{5mm}
$2$ & $\frac{3}{2}$ & $-2$ & $\frac{1}{2}$ & & & $90.00^\circ$\\[2mm]
$3$ & $\frac{11}{6}$ & $-3$ & $\frac{3}{2}$ & $-\frac{1}{3}$ & & $86.03^\circ$\\[2mm]
$4$ & $\frac{25}{12}$ & $-4$ & $3$ & $-\frac{4}{3}$ & $\frac{1}{4}$ & $73.35^\circ$\\[1mm]
\hline
\end{tabular}\\
\parbox{9cm}{
\caption{\small
Coefficients and stability angles for $L(\alpha)$-stability
of BDF methods for $s=2,3,4$.}
\label{tab-bdf-coeff}
}
\end{table}

\subsection{General IMEX-Peer Methods}
\subsubsection{The case $s=2$}
The two-stage singly implicit methods (\ref{implpeerq0}) with order two
form a one-parameter family, with free parameter $c_1$, say. The choice
$c_1=1/2$ produces the above implicit BDF2 method with step-size $\st/2$.
Note that requiring optimal zero stability, with $P$ having a single
eigenvalue one and the other zero, yields a completely defined method.
However, this would exclude interesting methods, such as the BDF2 method.

\begin{figure}[h]
\centering
\includegraphics[width=0.85\textwidth]{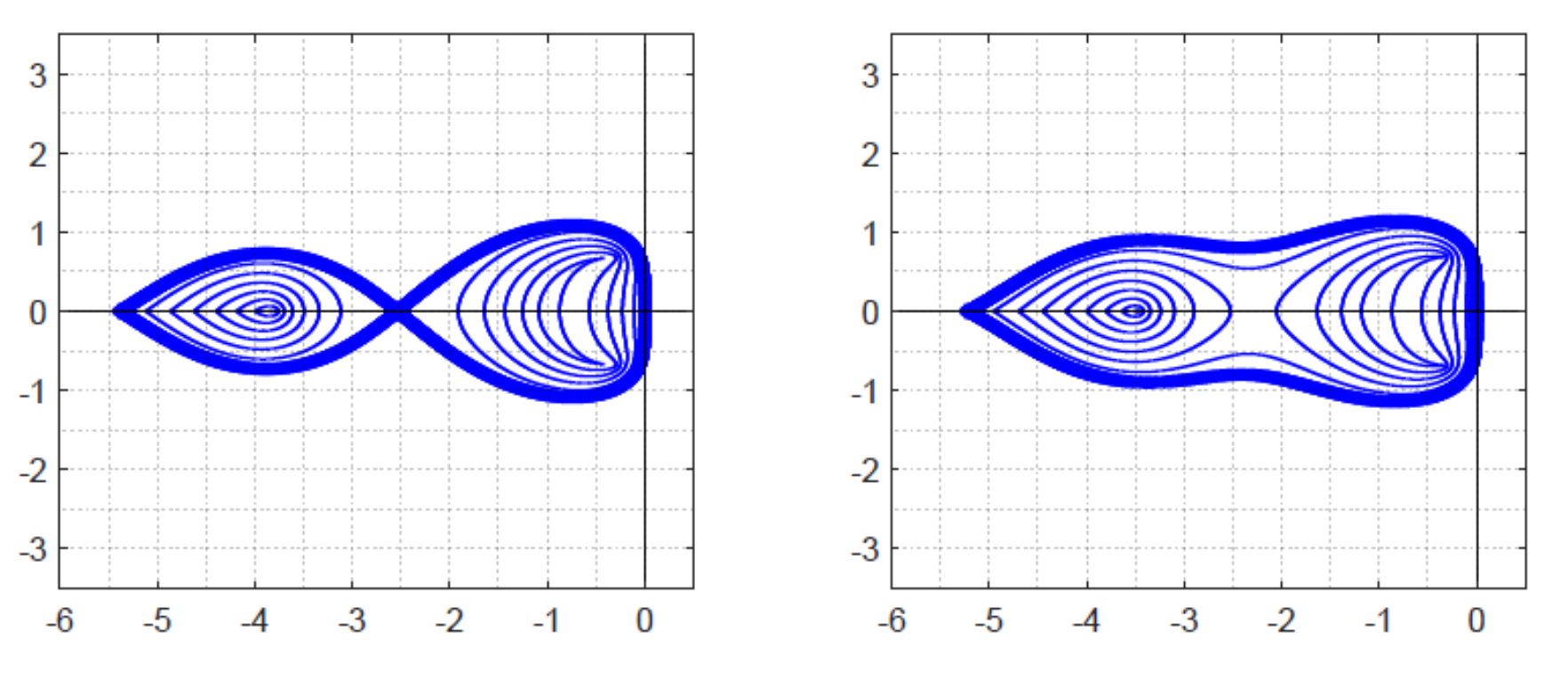}
\parbox{13cm}{
\caption{\small
IMEX-Peer2 method. Implicit method with $c_1=\frac12$, i.e., BDF2
with $\st/2$. Plots of the stability regions of the corresponding
explicit methods with $s_{21}^{(2)}=\mu\ne 0$; left panel with optimal
$\mu=\mu^*=10-4\sqrt{5}$, right panel with $\mu=\mu^*+\frac{1}{10}$.}
\label{fig-stab-peer2-exp}
}
\end{figure}

In order to find an IMEX method, where the explicit method has a larger
stability region, we start with the implicit BDF2 method with step size
$\st/2$ and then apply extrapolation with a strictly lower triangular
$S_2=(s_{ij}^{(2)})\ne 0$, say $s_{21}^{(2)}=\mu\ne 0$. Note that $\mu=2$
recovers the IMEX-BDF2 method from above. A careful study of the stability matrix
revealed that the largest interval $(-\beta_R,0)$ of the real negative
axis in the stability region is obtained if $\mu$ is the smallest root
of the polynomial $\mu^2-20\mu +20$, i.e.,
$\mu=\mu^*=10-4\sqrt{5}\approx 1.0557$ with real stability boundary
$\beta_R\approx 5.38$. Choosing $\mu$ equal to this optimal $\mu^*$ gives
a stability region which is pinched off at the real point
$x^*\approx -2.54$. Taking $\mu$ a bit larger, for example
$\mu=\mu^*+1/10$, gives a better shaped stability region, as shown in
Figure~\ref{fig-stab-peer2-exp}. The coefficients of the resulting IMEX-Peer2
method are
\begin{equation}
c =
\begin{pmatrix}
\frac{1}{2}\\[2mm]
1
\end{pmatrix},\quad
P =
\begin{pmatrix}
-\frac{1}{3} & \frac{4}{3}\\[2mm]
-\frac{4}{9} & \frac{13}{9}
\end{pmatrix},\quad
R =
\begin{pmatrix}
\frac{1}{3} & 0\\[2mm]
\frac{4}{9} & \frac{1}{3}
\end{pmatrix},\quad
S_2 =
\begin{pmatrix}
0 & 0\\[2mm]
\mu & 0
\end{pmatrix}\,,
\end{equation}
accomplished with $\hR=RS_2$ and $\hQ=R(I-S_2)V_0V_1^{-1}$.

\subsubsection{The cases $s=3$ and $s=4$}
To ensure good stability of the implicit Peer method with three or four
stages, we start with superconvergent singly-implicit methods
of order $p=s$ for constant step size and optimal zero-stability, developed
by Beck, Weiner, Podhaisky and Schmitt
(\cite{BeckWeinerPodhaiskySchmitt2012}, Table~3)
for large stiff ODE systems. These methods are $L(\alpha )$-stable with
angles $86.1^\circ$ and $83.2^\circ$, and possess relatively small error constants. The
nodes are selected such that $0<c_1<c_2<\ldots<c_s=1$. The
free parameters are then the $d=s(s-1)/2$ inputs of the matrix $S_2$.

As a design criterion, we would like to balance between optimal stability
regions $\St_\alpha$ and small error constants $c_{ex}$ for the extrapolation.
The latter one is very important since extrapolation has to be done forward in
time, i.e., future values are approximated outside the range of given time
points, which might cause relatively large errors. We expect that optimizing
$\St_\alpha$, i.e., maximising the size of its area, $|\St_\alpha|$, results
also in reasonably shaped stability regions $\St_E$ of the explicit methods.

Eventually, we perform an optimisation over the parameter space $p\in\R^d$ to compute
\begin{equation}
\label{obj_func}
p^* = \mbox{argmin} \{ -|\St_\alpha| + 1.5\,10^s\,|c_{ex}-c_0|\}\,,
\end{equation}
where $c_{ex}$ is defined in (\ref{c_ex}) and $c_0$ is the error constant
that corresponds to the extrapolation based on the $s$ most recently computed
stage values, i.e.,
\begin{equation}
\label{spec_nodes}
w_{n-1}^{(j+1)},\ldots,w_{n-1}^{(s)}=w_n^{(0)},w_n^{(1)},\ldots,w_n^{(j)},
\quad j=0,\ldots,s-1\,.
\end{equation}
In this case, $S_1$ is an upper triangular matrix and the relation $S_1V_1=(I-S_2)V_0$
is uniquely solvable for $S_1$ and $S_2$.
We find $c_0=4.1082\,10^{-2}, 4.2632\,10^{-3}$ for $s=3,4$, respectively.
The entries of the specific matrix $S_2$ are taken as initial guess
$p_0=(s^{(2)}_{21},s^{(2)}_{31},\ldots,s^{(2)}_{s,s-1})$
for the routine {\it fminsearch} implemented in {\sc Matlab} to compute the optimal $p^*$.
The choice of the objective function in (\ref{obj_func}) is motivated by two
requirements: (i) Due to the natural ordering of the nodes $c_i$, we take the nodes
given in (\ref{spec_nodes}) as reference set for a reasonable extrapolation
and aim at a moderate relative deviation of the error constants $c_{ex}$
in the range of nearly $10$ percent. (ii) The terms related to stability and
extrapolation have to be well balanced. This defines, with an anticipated target value
$|\St_\alpha|\approx 6$, which is slightly beat by the IMEX-BDF methods
(see Tab.~\ref{tab-stab-imex}), and the given
size of the constants $c_0$, the weighting factor $1.5\,10^s$ for the difference
$|c_{ex}-c_0|$.

The optimal parameters delivered by the minimisation process are
\begin{equation}
\begin{array}{lcllcl}
p_{21} &=& 4.6617853424698374\,10^0, & p_{31} &=& 3.3696230360366979\,10^0,\\
p_{32} &=& 5.6686050026329915\,10^{-1}, &&&
\end{array}
\end{equation}
for the IMEX-Peer3 method and
\begin{equation}
\begin{array}{lcrlcr}
p_{21} &=&  4.0913830614894255\,10^0, & p_{31} &=& -1.2244427616780204\,10^1,\\
p_{32} &=&  5.7564397758588521\,10^0, & p_{41} &=&  1.0587962913073733\,10^1,\\
p_{42} &=& -7.7409749651373776\,10^0, & p_{43} &=&  4.1019377658951353\,10^0,
\end{array}
\end{equation}
for the IMEX-Peer4 method. We set $s_{ij}^{(2)}=p_{ij}$ to define $S_2$ in each case.

The resulting values for the stability regions $\St_\alpha$ and $\St_E$ as well
as for the error constants are collected in Table~\ref{tab-stab-imex}. For
comparison, we also show the values for the IMEX-BDF methods. It can be
observed that (i) the error constants for the extrapolation are comparable,
(ii) the sizes of the stability regions differ only moderately, and (iii) the
IMEX-Peer methods have a significantly larger interval (up to a factor two
for the two-stage method) on the negative real axis included in the stability
region. More details are visible in Figure~\ref{fig-stabreg-all-shape}.

\begin{table}[h]
\centering
\begin{tabular}{|l|crr|rr|cc|}
\hline\rule{0mm}{5mm}\hspace{-0.1cm}
Method & $\alpha$ & $|\St_\alpha|$ & $x_{max}$ & $|\St_E|$ & $x_{max}$ & $c_{im}$ & $c_{ex}$\\[1mm]
\hline\rule{0mm}{5mm}\hspace{-0.1cm}
IMEX-BDF2 & $90.0^\circ$ & $6.28$ & $-2.67$ & $6.98$ & $-2.67$ & $7.05\,10^{-2}$ & $2.11\,10^{-1}$\\[2mm]
IMEX-BDF3 & $86.0^\circ$ & $7.27$ & $-2.86$ & $9.65$ & $-2.86$ & $8.93\,10^{-3}$ & $3.57\,10^{-2}$\\[2mm]
IMEX-BDF4 & $73.4^\circ$ & $7.30$ & $-2.84$ & $9.92$ & $-2.84$ & $8.91\,10^{-4}$ & $4.45\,10^{-3}$\\[1mm]
\hline\rule{0mm}{5mm}\hspace{-0.1cm}
IMEX-Peer2  & $90.0^\circ$ &  $7.44$ & $-4.86$ &  $8.53$ & $-5.22$ & $7.05\,10^{-2}$ & $2.78\,10^{-1}$\\[2mm]
IMEX-Peer3  & $86.1^\circ$ &  $8.28$ & $-3.07$ & $10.68$ & $-3.07$ & $8.20\,10^{-3}$ & $3.58\,10^{-2}$\\[2mm]
IMEX-Peer4  & $83.2^\circ$ &  $4.64$ & $-3.57$ & $9.36$ & $-3.57$ & $3.43\,10^{-4}$ & $4.27\,10^{-3}$\\[1mm]
\hline
\end{tabular}\\
\parbox{13cm}{
\caption{\small
Size of stability regions $\St_\alpha$ and $\St_E$, $x_{max}$ at the negative
real axis and error constants for IMEX-BDF and IMEX-Peer methods with $s=2,3,4$.}
\label{tab-stab-imex}
}
\end{table}

\section{Comparison of Stability Regions}
\begin{figure}[h]
\setlength{\unitlength}{1cm}
\centering
\includegraphics[width=0.32\textwidth]{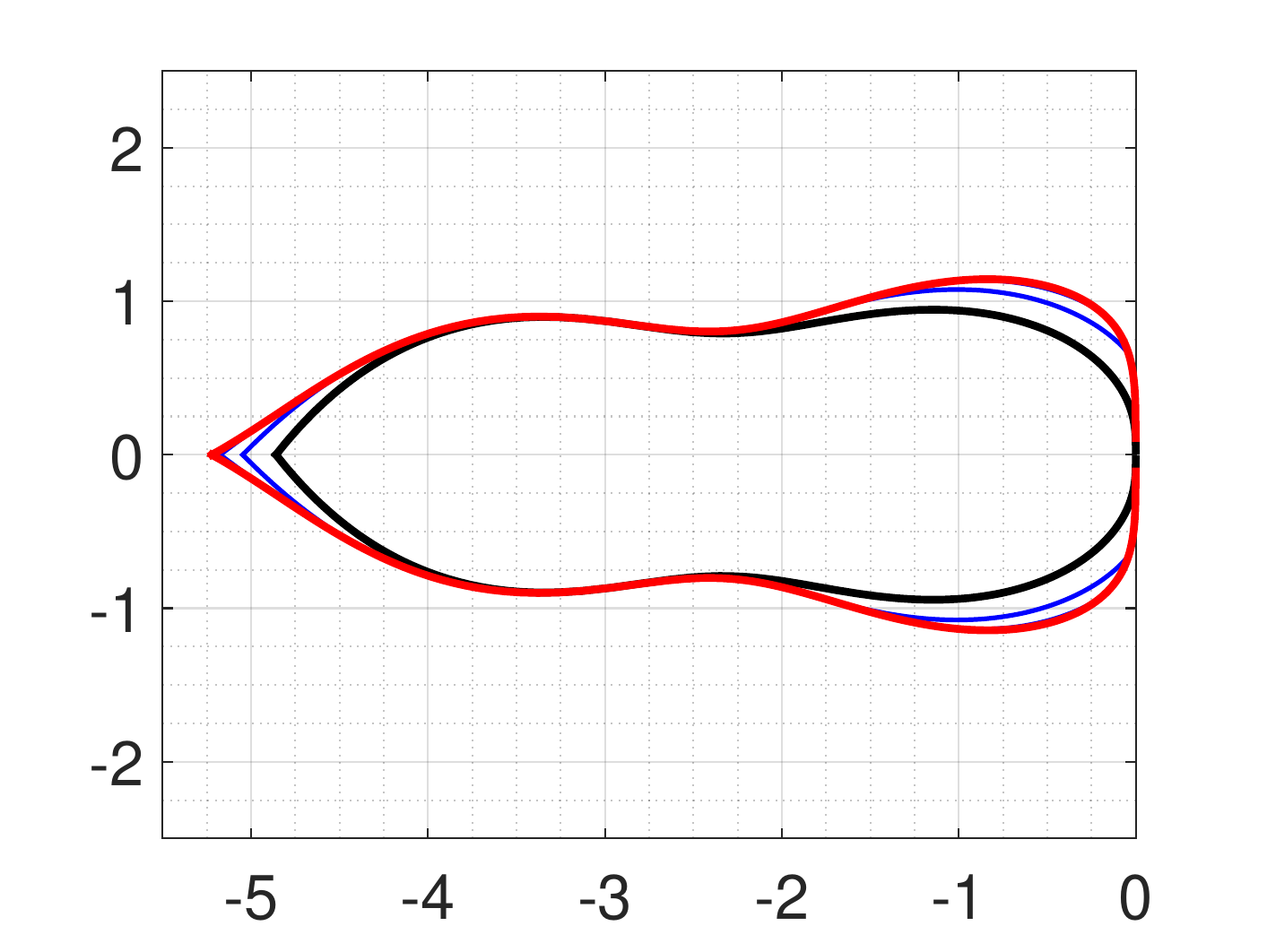}
\put(-3.8,2.8){\scriptsize IMEX-PEER2}
\includegraphics[width=0.32\textwidth]{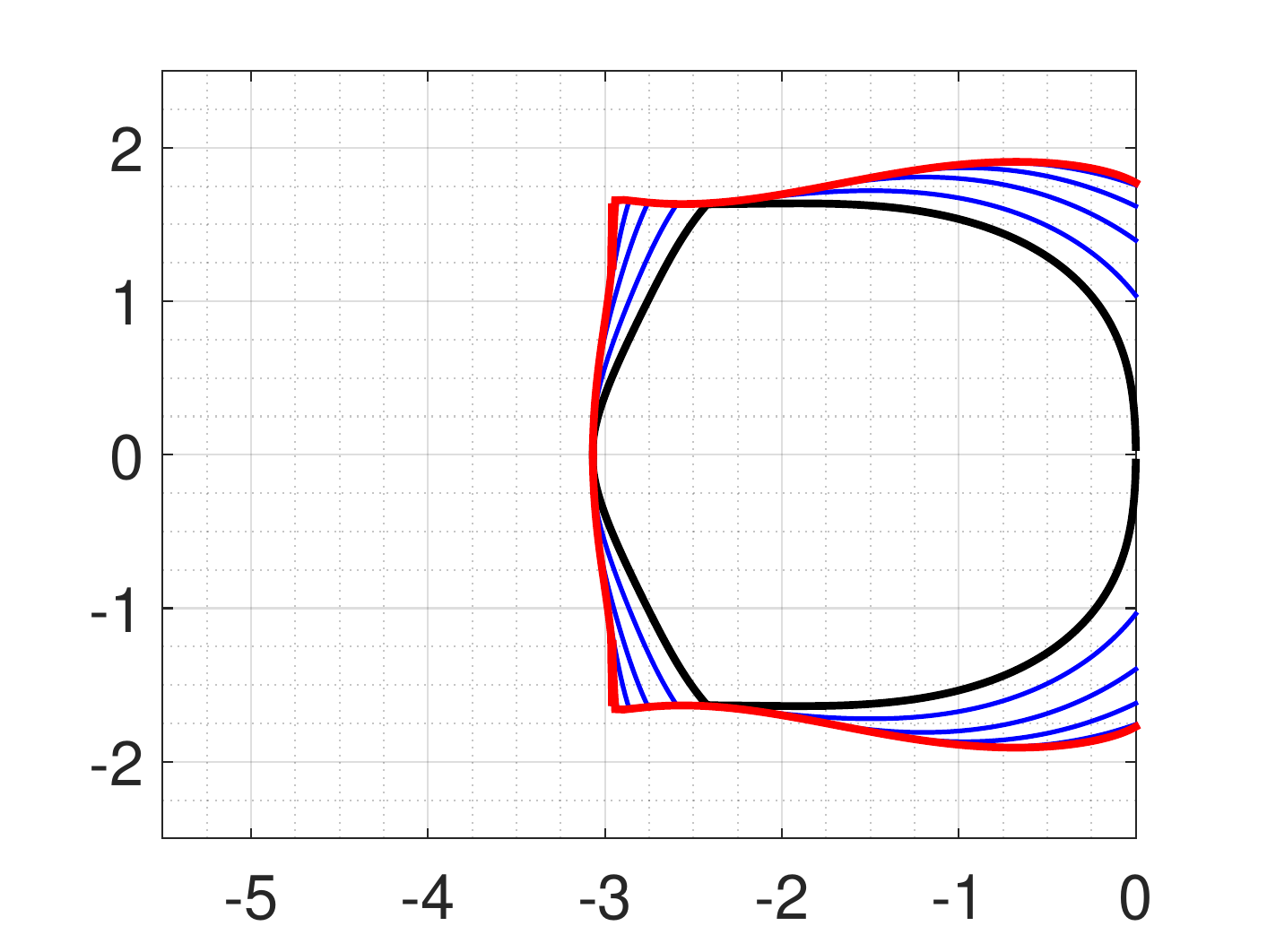}
\put(-3.8,2.8){\scriptsize IMEX-PEER3}
\includegraphics[width=0.32\textwidth]{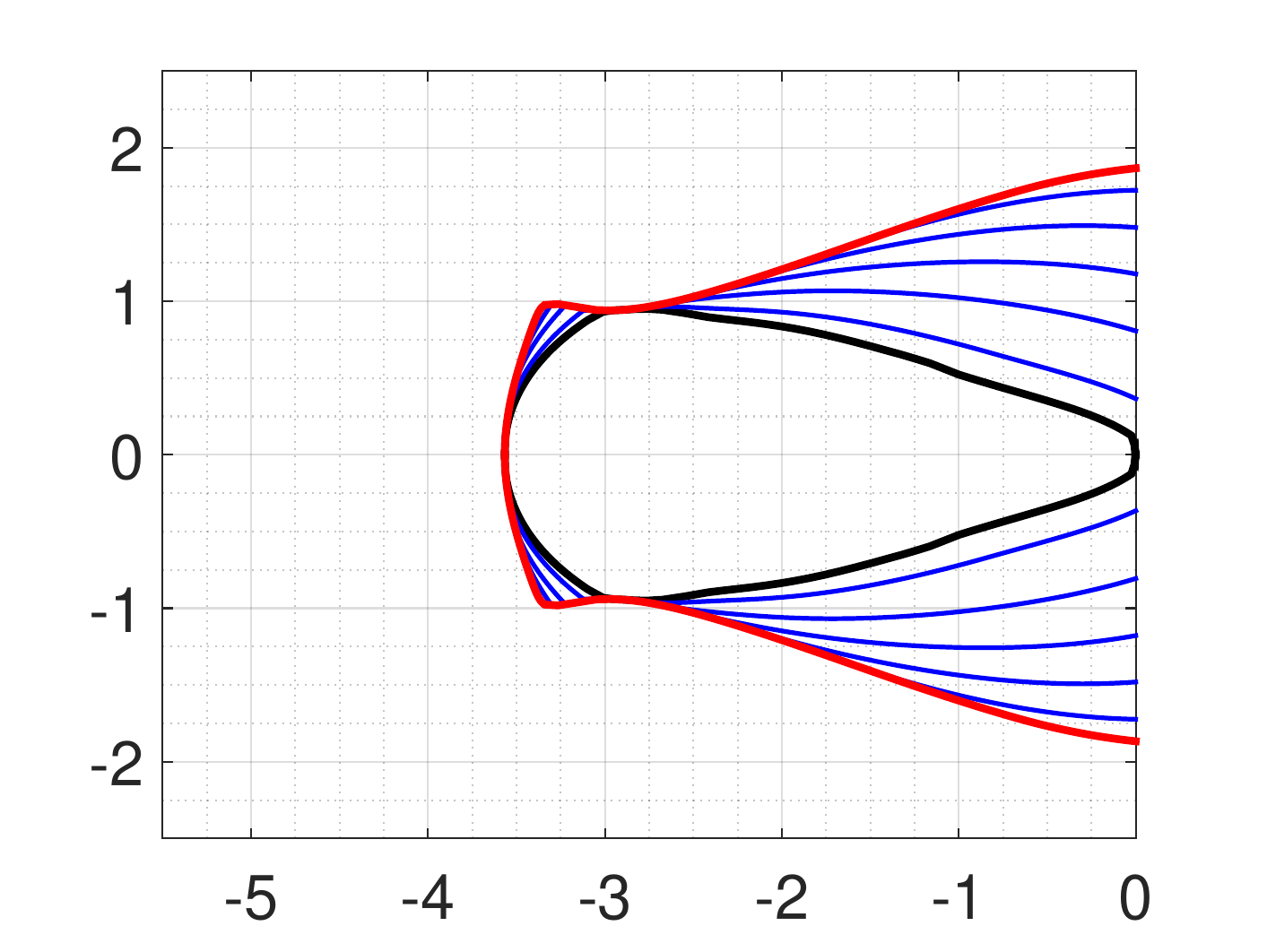}
\put(-3.8,2.8){\scriptsize IMEX-PEER4}
\\
\includegraphics[width=0.32\textwidth]{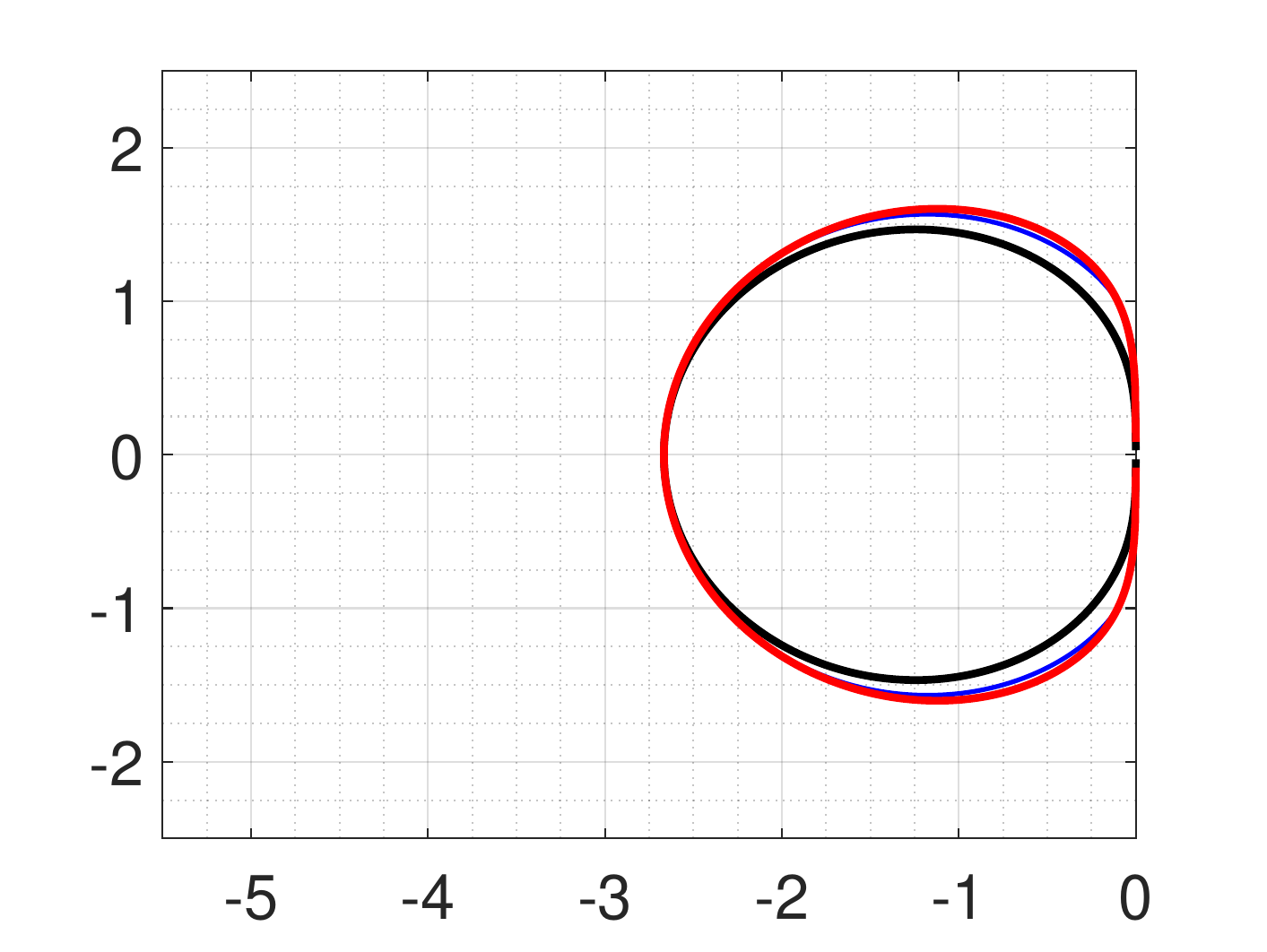}
\put(-3.8,2.8){\scriptsize IMEX-BDF2}
\includegraphics[width=0.32\textwidth]{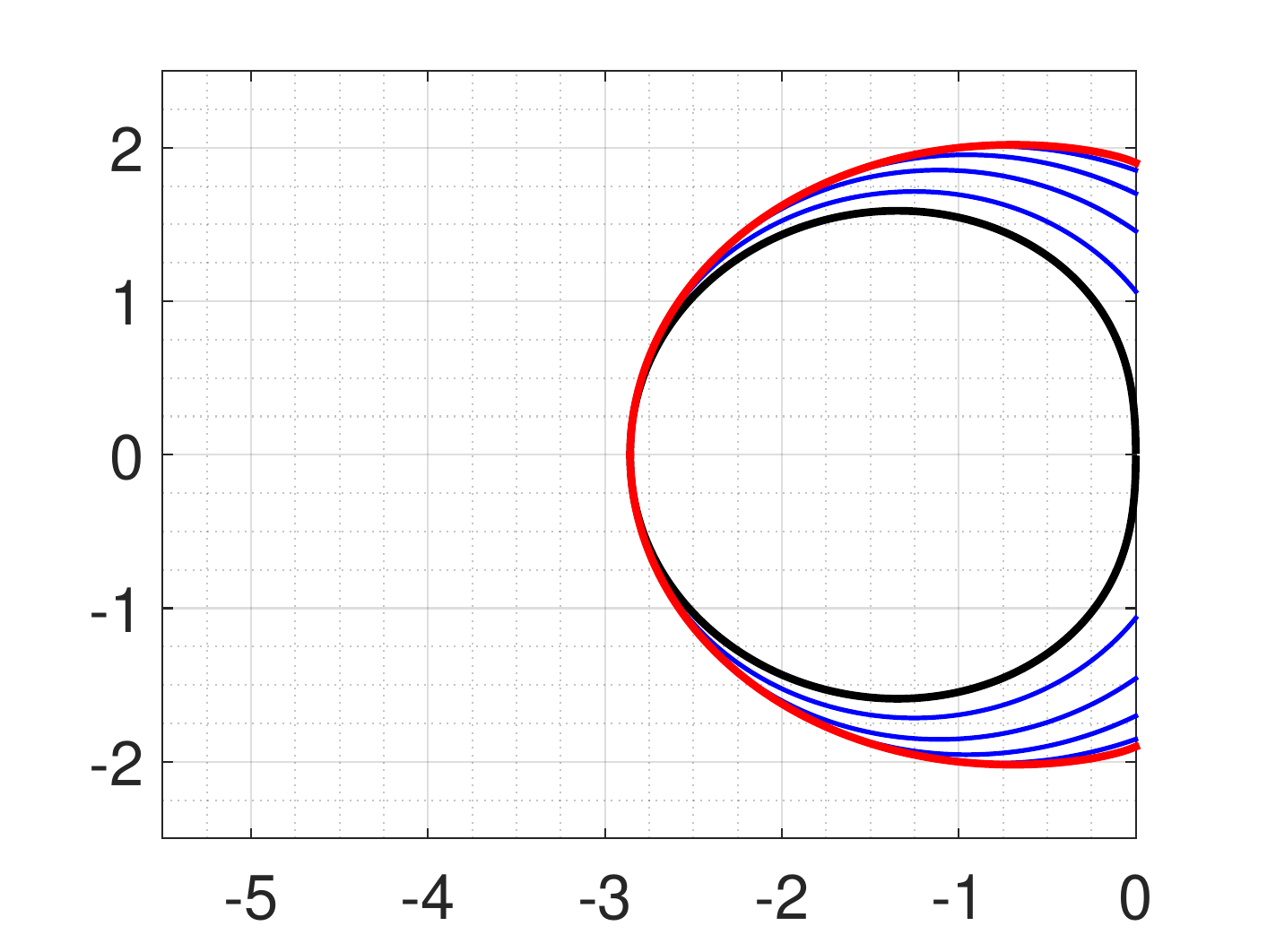}
\put(-3.8,2.8){\scriptsize IMEX-BDF3}
\includegraphics[width=0.32\textwidth]{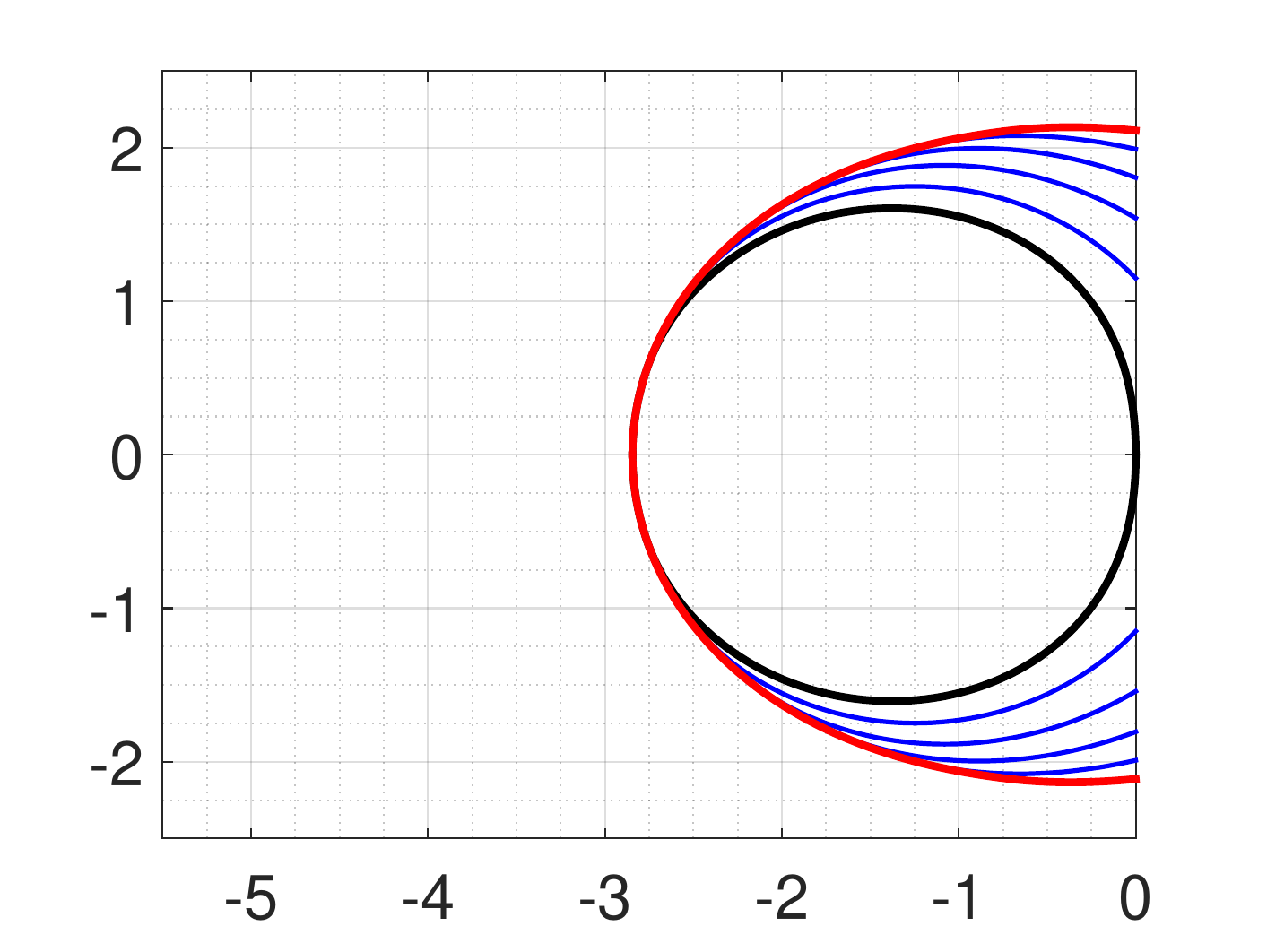}
\put(-3.8,2.8){\scriptsize IMEX-BDF4}
\\
\includegraphics[width=0.32\textwidth]{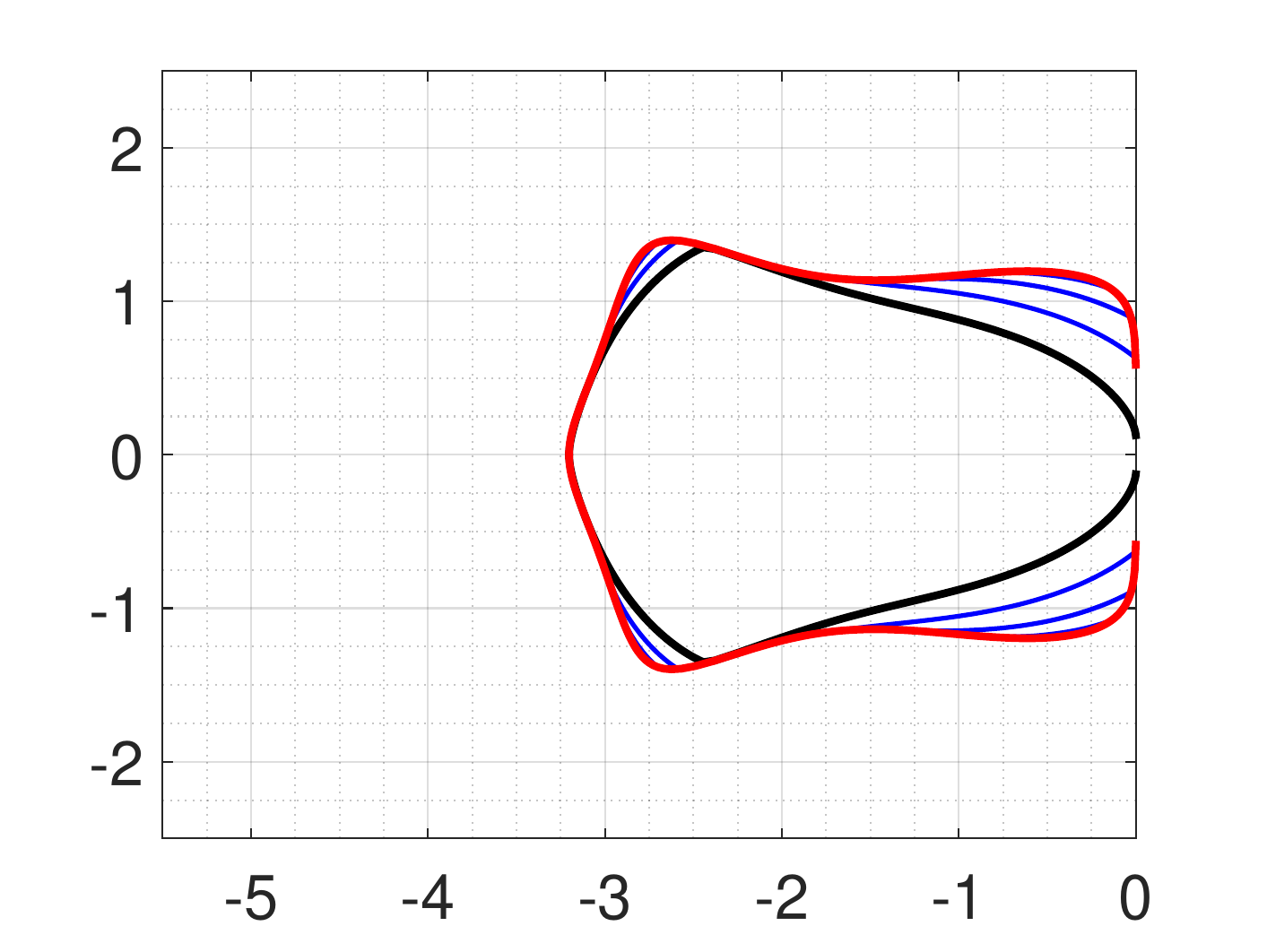}
\put(-3.8,2.8){\scriptsize IMEX-DIMSIM2}
\includegraphics[width=0.32\textwidth]{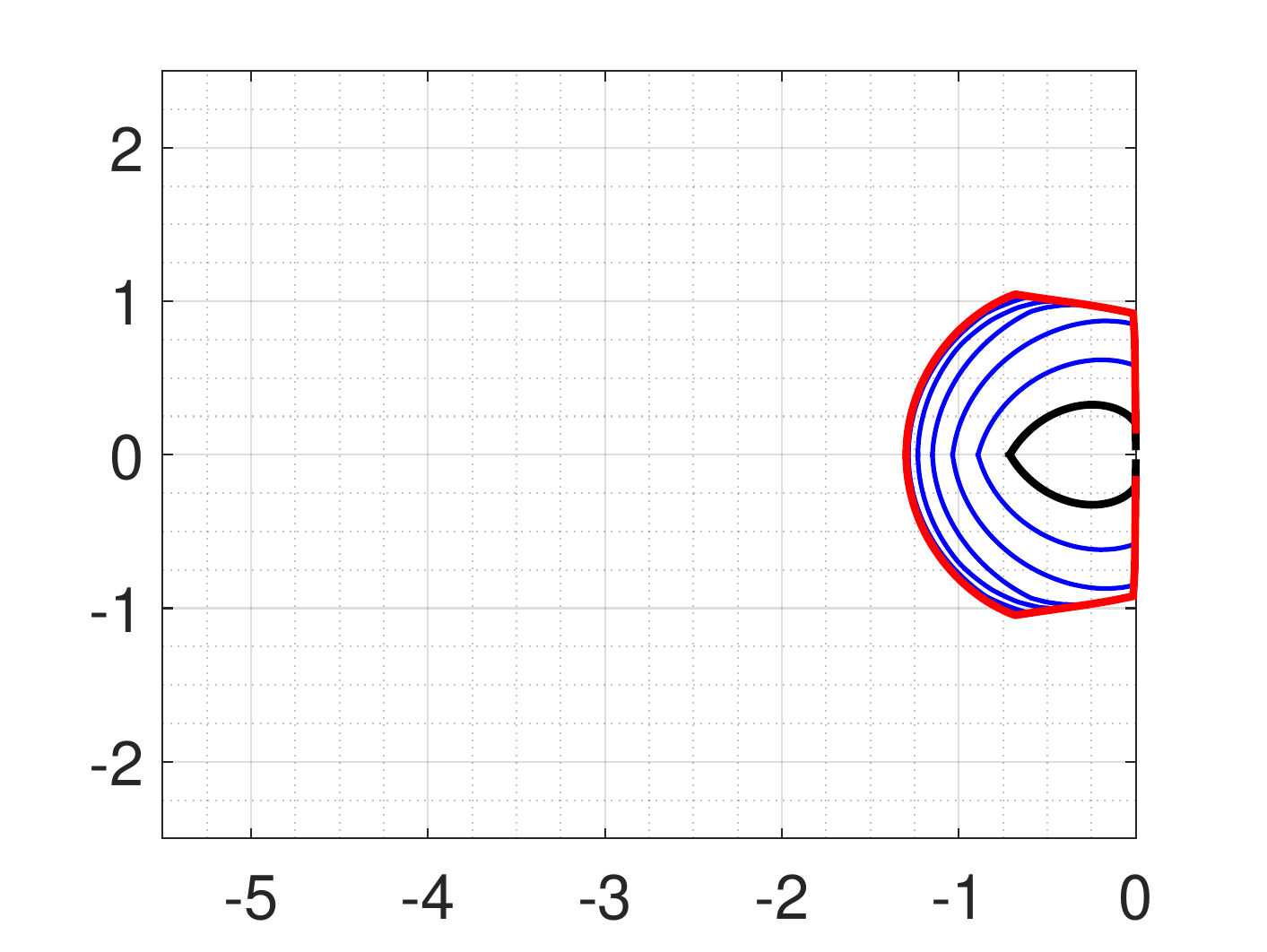}
\put(-3.8,2.8){\scriptsize IMEX-DIMSIM3}
\includegraphics[width=0.32\textwidth]{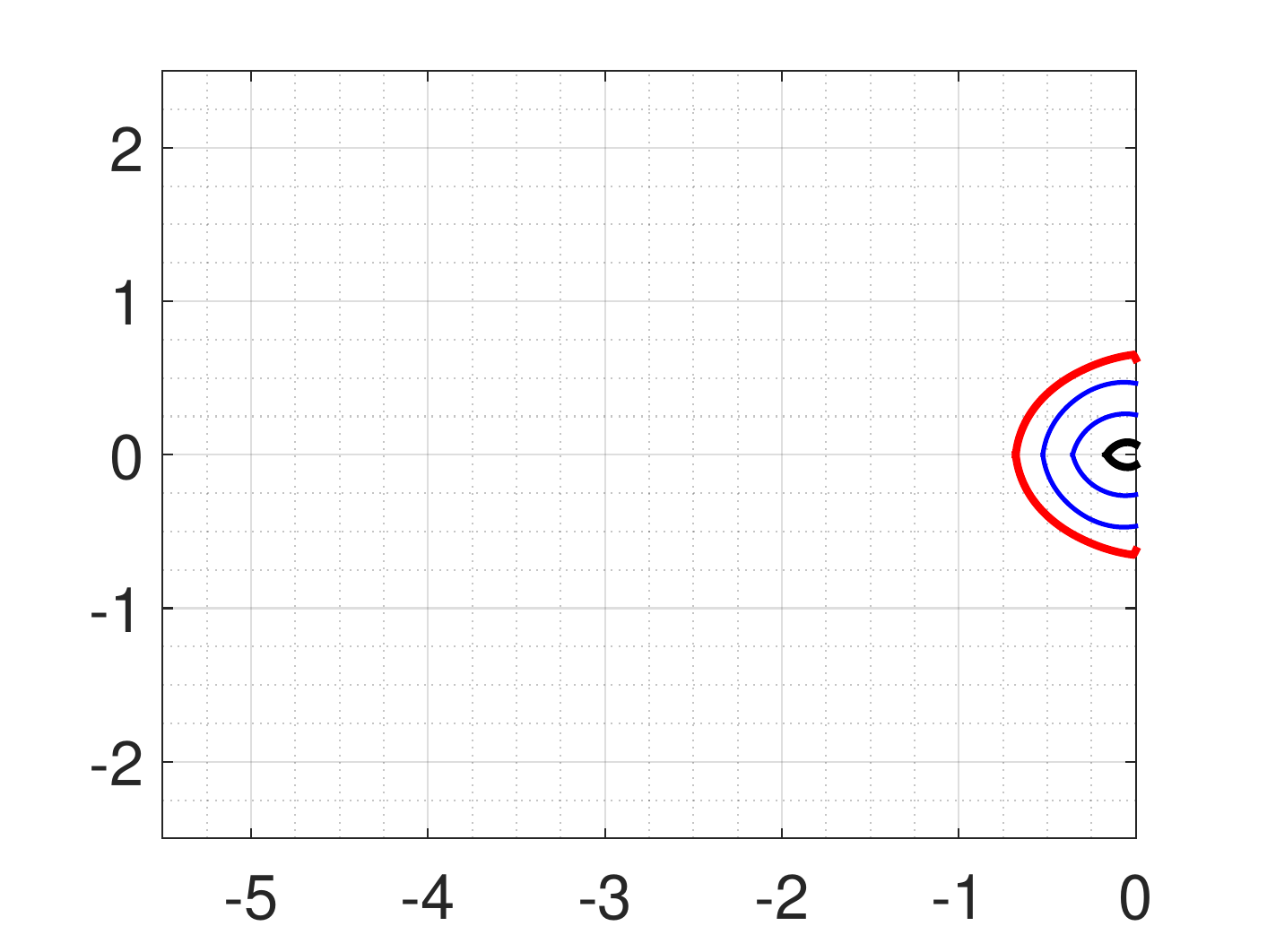}
\put(-3.8,2.8){\scriptsize IMEX-DIMSIM4}
\\
\parbox{13cm}{
\caption{
Top to bottom: stability regions
$\St_{\alpha}$ (black line), $\St_{\beta}$ for
$\beta=75^\circ,60^\circ,45^\circ,30^\circ,15^\circ$ (blue lines),
and $\St_0$ (red line) for IMEX-PEER(s), IMEX-BDF(s), and IMEX-DIMSIM(s)
methods, $s=2,3,4$ (left to right).}
\label{fig-stabreg-all-shape}
}
\end{figure}

Here we compare the stability regions of the IMEX-Peer and IMEX-BDF methods to those
of the IMEX-DIMSIM methods developed and tested by Cardone et al. \cite{CardoneJackiewiczSanduZhang2014}.
There, the authors first selected an implicit DIMSIM method with suitable stability and order
properties, and then obtained the explicit component through an optimisation procedure
that maximized the combined stability region of the pair. They applied this methodology
to construct IMEX pairs of orders one to four. In contrast, we took also care of the
error constants for the underlying extrapolation process. The stability regions $\St_\beta$ for
varying angle $\beta$ and methods with $s=2,3,4,$ are shown in Figure~\ref{fig-stabreg-all-shape}.

\begin{figure}[h]
\centering
\includegraphics[width=0.9\textwidth]{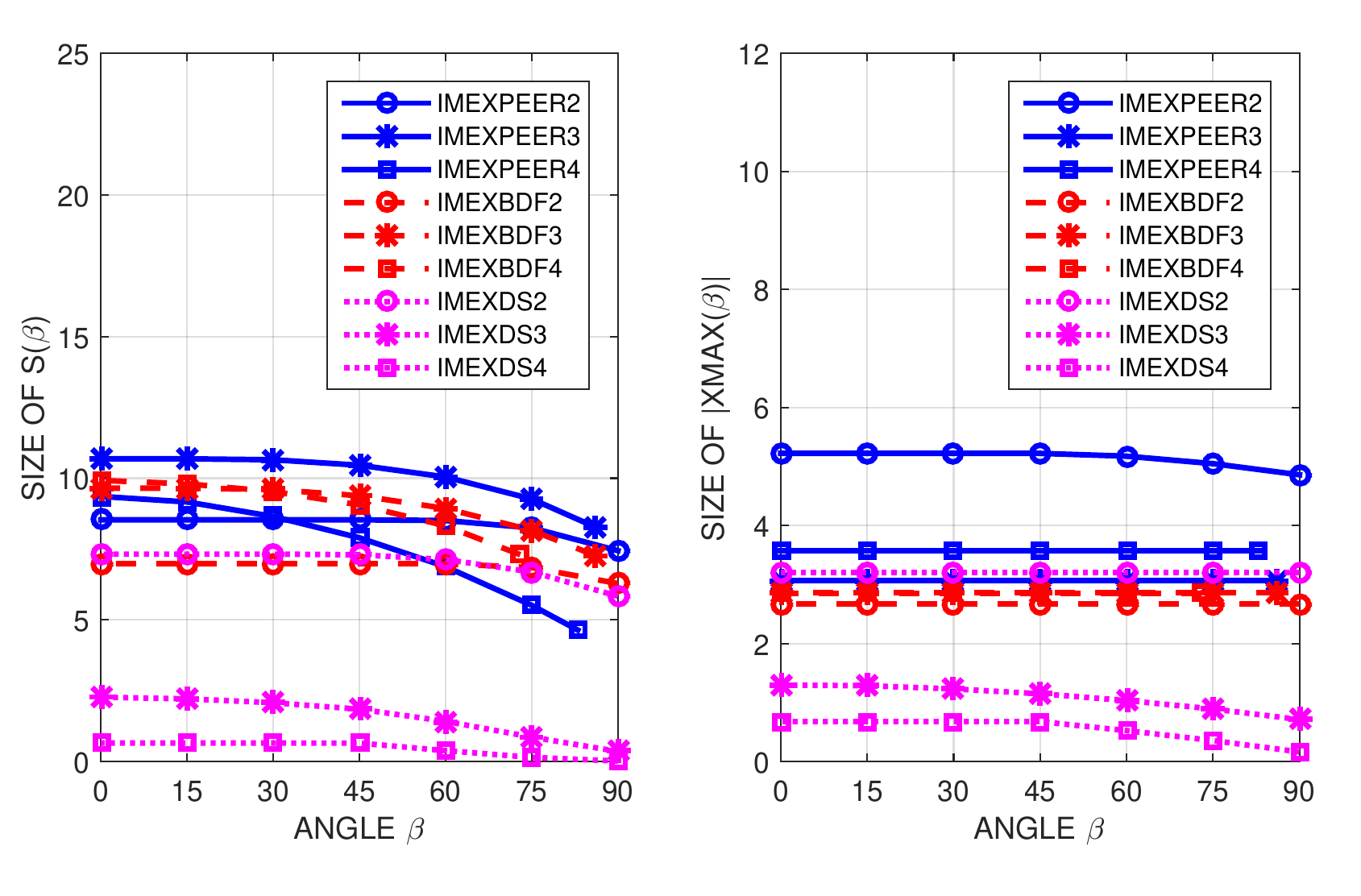}
\parbox{13cm}{
\caption{\small
Size of stability regions $\St_{\beta}$ for
$\beta=\alpha,75^\circ,60^\circ,45^\circ,30^\circ,15^\circ,0^\circ$
(left panel) and $|x_{max}|$ on the negative real axis (right panel) for IMEX-PEER(s),
IMEX-BDF(s), and IMEX-DIMSIM(s) methods, $s=2,3,4$.}
\label{fig-stabreg-all-size}
}
\end{figure}

It is obvious that the two-step methods of Peer type allow the construction of higher-order
extrapolation-based IMEX schemes with larger stability regions. Whereas the IMEX-DIMSIM2 scheme
is still competitive with respect to absolute size, the other two IMEX-DIMSIM schemes suffer
clearly from small stability regions. For these methods, we expect stability problems for
larger time steps, which is indeed confirmed by our numerical experiments.
In Figure~\ref{fig-stabreg-all-size}, we have collected the values for the size of
stability regions $\St_\beta$ with $\beta=\alpha,75^\circ,60^\circ,45^\circ,30^\circ,15^\circ$ and the absolute
value of the left-most point $x_{max}\in \St_E$ on the negative real axis.

\section{Numerical Experiments}

\subsection{Linear Advection-Reaction Problem}
A first PDE problem for accuracy test is a linear advection-reaction system from
\cite{HundsdorferRuuth2007}. The equations are
\begin{eqnarray}
\partial_t u + \alpha_1\,\partial_x u &=& -k_1 u + k_2 v + s_1\,,\\
\partial_t v + \alpha_2\,\partial_x v &=& k_1 u - k_2 v + s_2
\end{eqnarray}
for $0<x<1$ and $0<t\le 1$, with parameters
\[
\alpha_1=1,\;\alpha_2=0,\;k_1=10^6,\;k_2=2k_1,\;s_1=0,\;s_2=1,
\]
and with the following initial and boundary conditions:
\[
u(x,0)=1+s_2x,\;v(x,0)=\frac{k_1}{k_2}u(x,0)+\frac{1}{k_2}s_2,\;
u(0,t)=1-\sin(12t)^4\,.
\]
Note that there are no boundary conditions for $v$ since $\alpha_2$
is set to be zero.

Fourth-order finite differences on a uniform mesh consisting of $m=400$
nodes are applied in the interior of the domain. At the boundary we can
take third-order upwind biased finite differences, which here does not
affect an overall accuracy of four \cite{HundsdorferRuuth2007}
and gives rise for a spatial error of $1.5\,10^{-5}$.

\begin{figure}[ht]
\centering
\includegraphics[width=0.8\textwidth]{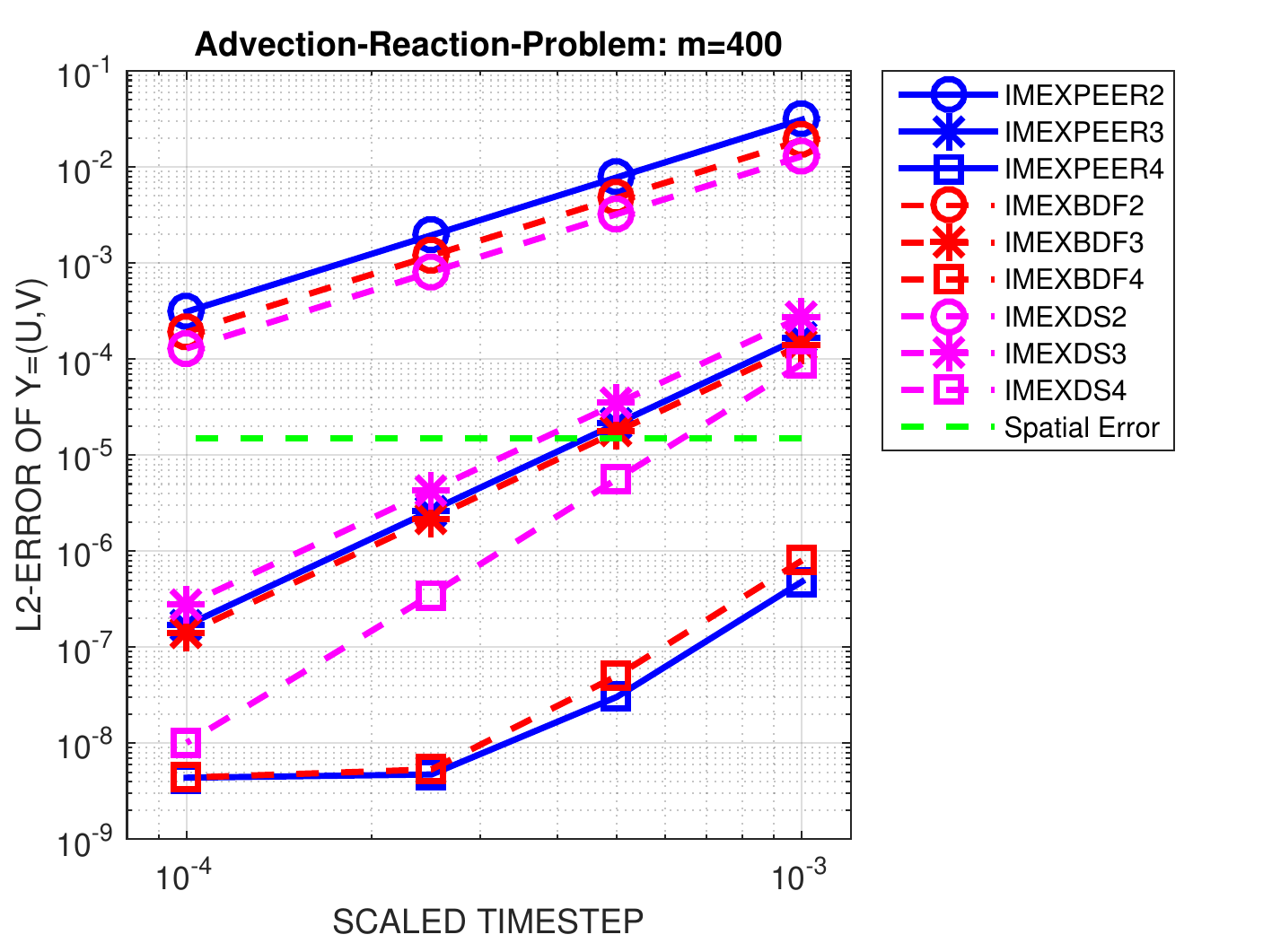}
\parbox{13cm}{
\caption{\small
Advection-Reaction-Problem: Temporal $l_2$-errors at $T=1$
of the total concentration vs. scaled step sizes, $m=400$.
Comparison of IMEX-Peer, IMEX-BDF and IMEX-DIMSIM methods.}
\label{fig-advreac400}
}
\end{figure}

In the IMEX setting, the reaction is treated implicitly and all other
terms explicitly. Accurate initial values are computed by the variable
step-size code ODE15S with high tolerances. We have used step sizes
$\st=10^{-3},\,5\,10^{-4},\,2.5\,10^{-4},\,10^{-4}$ and compared the
numerical values at the final time $T=1$ with an accurate reference
solution in the $l_2$-vector norm as in \cite{CardoneJackiewiczSanduZhang2014}.
The results are plotted in Figure~\ref{fig-advreac400}.

All second-order and third-order methods show their classical
orders and perform nearly identical for this problem. For the
IMEX-DIMSIM4 method we observe order four, but the larger error
constant compared to the IMEX-Peer4 and IMEX-BDF4 scheme is apparent.
The similar asymptotic behaviour for the latter shows an order reduction,
which was also observed in \cite{HundsdorferRuuth2007} as an inherent
issue for very high-accuracy computations. However, this effect
appears on a level far below the spatial discretisation error.

\subsection{Nonlinear Adsorption-Desorption Problem}
The problem is taken from \cite{HundsdorferRuuth2007}. Let $u$ and $v$
be the dissolved and adsorbed concentration, respectively, satisfying
the equations
\begin{eqnarray}
\partial_t u + a\,\partial_x u &=& \hspace{0.3cm}\kappa (v - \phi(u))\,,\\
\partial_t v &=& -\kappa (v - \phi(u))
\end{eqnarray}
for $0<x<1$ and $0<t\le T$, with $\phi(u)=k_1u/(1+k_2u)$. The initial values are
set to zero, $u_0=v_0=0$, and an oscillatory inflow condition is taken to get some
smooth variations in the solution, along with the shocks:
\[
\begin{array}{rll}
u(0,t) &=1 - \cos^2(6\pi t) & \mbox{ if } a>0 \,,\\
u(1,t) &=0 & \mbox{ if } a<0\,.
\end{array}
\]
The parameters are $\kappa=10^6$, $k_1=50$, $k_2=100$, $T=1.25$, and the velocity is set to
\[ a = -\frac{3}{\pi}\arctan(100(t-1))\,, \]
giving approximately $a=1.5$ for $t<1$ (adsorption phase) and $a=-1.5$ for
$t>1$ (desorption phase).

We use the WENO5 scheme for the spatial discretisation from Shu (\cite{Shu1999}, formulas
$(2.58)-(2.63)$ with parameter $\varepsilon=10^{-12}$) on a uniform (cell centred) grid,
$x_i=(i-\frac12)\triangle x,\,i=1,\ldots,m$, with mesh width $\triangle x=1/m$.
This WENO5 spatial scheme provides high accuracy in smooth regions together
with good monotonicity properties near shocks. We set $m=800$ and note that in this
case the spatial error is $1.2\,10^{-3}$.

\begin{figure}[h]
\centering
\includegraphics[width=0.8\textwidth]{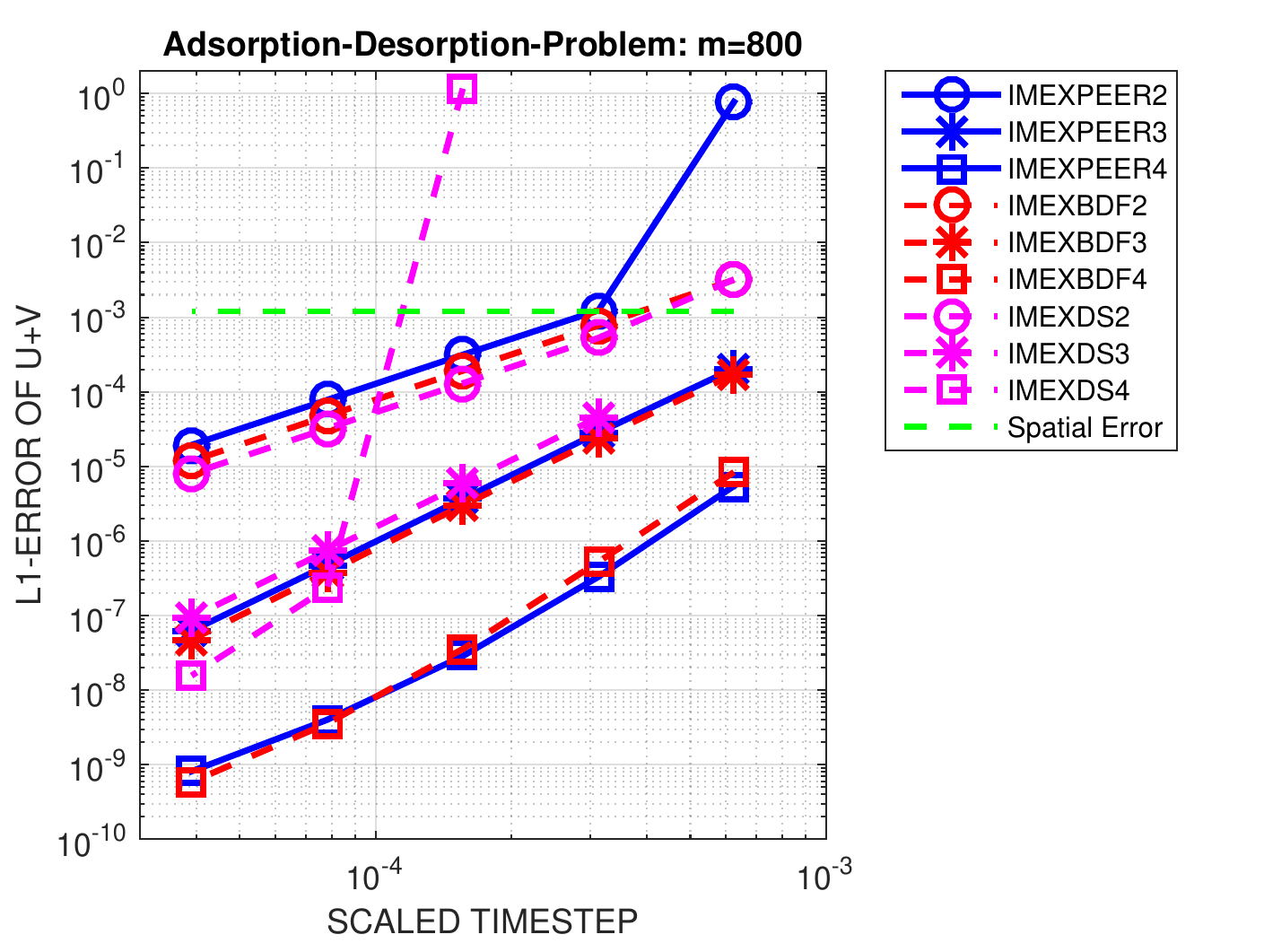}
\parbox{13cm}{
\caption{\small
Adsorption-Desorption-Problem: Temporal discrete $l_1$-errors of the total concentration
vs. scaled step sizes, $m=800$.
Comparison of IMEX-Peer, IMEX-BDF ans IMEX-DIMSIM methods.}
\label{fig-adsdes800}
}
\end{figure}

In the IMEX methods, the advection term is treated explicitly and the stiff relaxation
term implicitly, where a Newton method is efficiently performed at each spatial node separately.
The starting values for the methods are taken as $w_0=0$.
To allow a direct comparison to numerical schemes presented in \cite{HundsdorferRuuth2007},
we have used step sizes $\triangle\,t=2^{-j}\triangle\,x,j=1,\ldots,5,$ and compared the
numerical values of the total concentration, $u+v$, at the final time $t=T$ with an accurate
reference solution in the discrete $l_1$-norm ($\|v\|_1=h\sum_i|v_i|$),
see Figure~\ref{fig-adsdes800}.

As before, the results for the IMEX schemes with $s=2,3,$ largely coincide. We note
that the IMEX-BDF2 and IMEX-DIMSIM3 method did not converge for the largest time step.
We clearly observe stability problems for IMEX-DIMSIM4, which can be explained by
the relatively small stability region of the underlying explicit methods. The method needs
small time steps to prevent instabilities, and even then the error behaviour favours
the other fourth-order methods. IMEX-Peer4 and IMEX-BDF4 gave nearly identical results
with an increasing order reduction which was already visible in the first test problem.
In view of the spatial error, temporal errors below $10^{-4}$ are of less importance
for the total PDE error, however.

\subsection{The Schnakenberg Problem}
A classical example of two-dimensional reaction-diffusion equations for testing
numerical algorithms is the Schnakenberg system \cite{Schnakenberg1979,HundsdorferVerwer2003}.
The equations read
\begin{eqnarray}
\partial_t u &=& D_1\nabla^2 u + \kappa (a - u + u^2v)\,,\\
\partial_t v &=& D_2\nabla^2 v + \kappa (b - u^2v)\,,
\end{eqnarray}
where $u$ and $v$ denote the concentration of activator and inhibitor, respectively.
We follow the setup in \cite{HundsdorferVerwer2003} and take $D_1=0.05$, $D_2=1$, $\kappa=100$,
$a=0.1305$, $b=0.7695$, $T=1$. The solution is computed on the unit square domain
$\Omega=(0,1)^2$ with the initial conditions
\[
\begin{array}{rl}
u(x,y,0) &=  a+b+10^{-3}\,\exp\left(-100*((x-\frac13)^2+(y-\frac12)^2)\right)\,,\\[2mm]
v(x,y,0) &= \displaystyle\frac{b}{(a+b)^2}\,.
\end{array}
\]
and homogeneous Neumann boundary conditions.

For the spatial discretisation, we apply second-order finite
differences on a uniform (cell centred) grid,
$\left(x_i=(i-\frac12)h,y_j=(i-\frac12)h\right)$, $i,j=1,\ldots,m$, with mesh width
$h=1/m$, where $m=400$ has been taken.

\begin{figure}[h]
\centering
\includegraphics[width=0.8\textwidth]{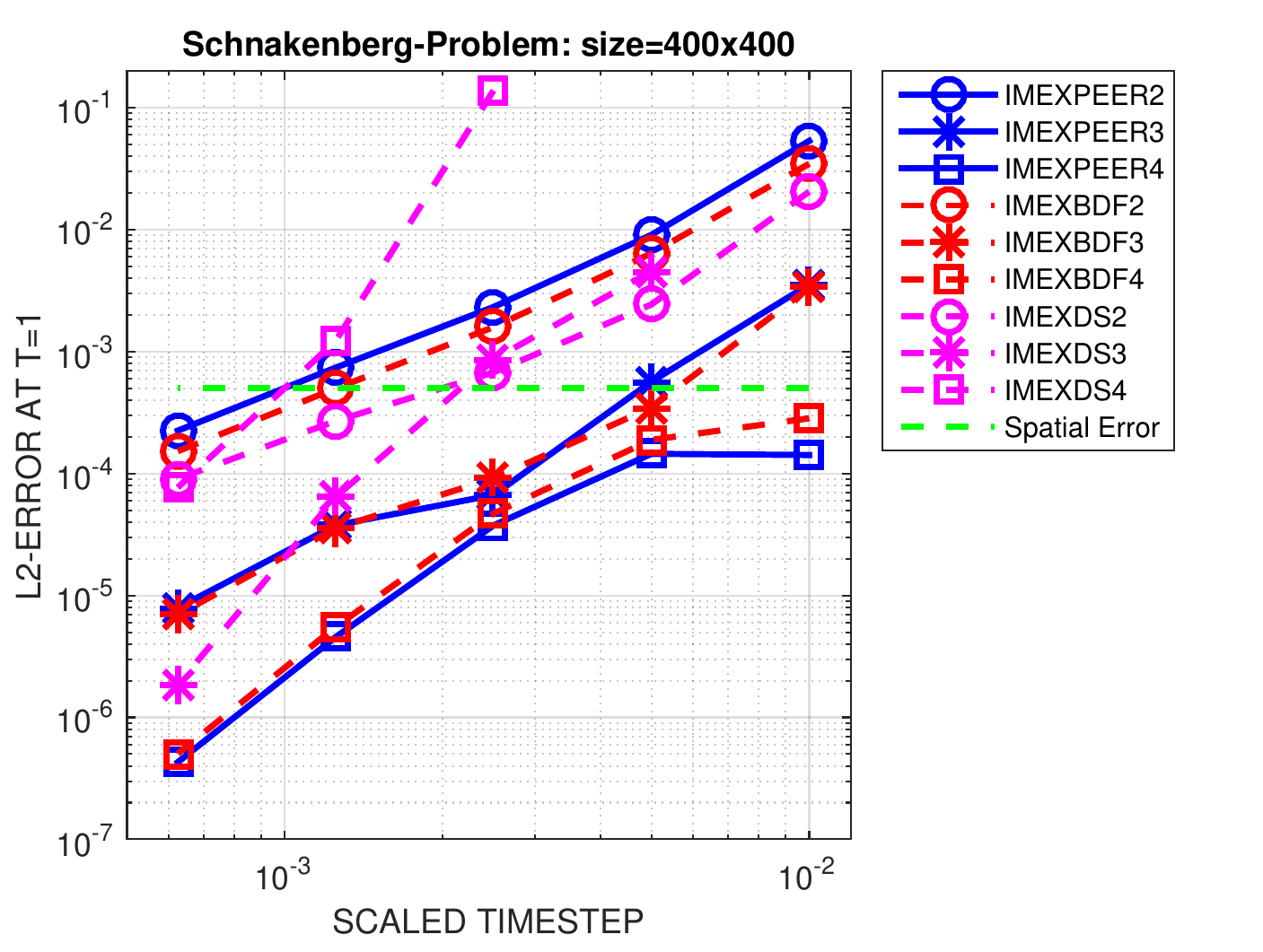}
\parbox{13cm}{
\caption{\small
Schnakenberg-Problem: Temporal discrete $l_2$-errors at $T=1$
vs. scaled step sizes, $m=400$.
Comparison of IMEX-Peer, IMEX-BDF and IMEX-DIMSIM methods.}
\label{fig-reacdiff400}
}
\end{figure}

Here, we treat the reaction explicitly and the diffusion implicitly. Accurate initial
values are computed by the variable step-size code ODE15S with high tolerances. We have
used step sizes $\triangle\,t=2^{3-j}/m,j=1,\ldots,5,$ and compared the numerical values
at the final time $t=T$ with an accurate reference solution in the discrete $l_2$-norm,
($\|v\|_2=\sqrt{h\sum_i|v_i|^2}$), see Figure~\ref{fig-reacdiff400}.

The second-order IMEX methods perform well and show their
classical order. IMEX-DIMSIM2 produced the best results due to
a smaller error constant. The higher-order IMEX-DIMSIM schemes
failed for larger time steps, whereas IMEX-DIMSIM4 gave again
unsatisfactory results at all. A further time, Peer and BDF
methods delivered nearly identical numerical solutions. Both
showed a somehow unpredictable behaviour for larger time steps,
but in this case they are still more efficient than the DIMSIM
schemes.

\section{Conclusion}
We have developed a new family of $s$-stage implicit-explicit
Peer methods, starting with $L(\alpha)$-stable implicit Peer
methods with order $p=s$ and applying an extrapolation of the
same order to preserve the order of convergence. The well-known
IMEX-BDF($s$) methods applied with constant step size $\st/s$ fit
into this framework when they are considered as $s$-stage Peer
methods with equidistant nodes and step size $\st$. We gave the
corresponding formula to convert. We examined the linear stability
properties of these IMEX methods to construct new IMEX-Peer
methods of order $p=2,3,4$, with optimally balanced size of
stability regions and error constants for the underlying extrapolation.
A detailed comparison with the recently proposed IMEX-DIMSIM methods
\cite{CardoneJackiewiczSanduZhang2014} showed a significant improvement
of the stability properties and a better performance of the
higher-order methods for three numerical test problems.

We are planning to extend this work to a variable step size environment
and to include other classes of implicit Peer methods, e.g., those
with strictly diagonal matrix $R$ to allow an efficient parallelisation.
We will also consider linearly implicit Peer methods of higher order
$p\ge 4$ as developed in \cite{PodhaiskyWeinerSchmitt2005} and
successfully applied in \cite{GerischLangPodhaiskyWeiner2009} to large scale
PDE problems within an adaptive Rothe approach, i.e., first discretise
in time and then apply an adaptive spatial discretisation afterwards.
There are L($\alpha$)-stable
methods of this type with reasonable large angles $\alpha$ and small error
constants available up to order $p=8$, which give them a clear advantage over
higher order BDF methods.

\section{Acknowledgement}
The authors would like to thank A.~Sandu and H.~Zhang for making
the precise parameters for the implicit-explicit DIMSIM methods
developed and tested in \cite{CardoneJackiewiczSanduZhang2014}
available for our comparison. Jens Lang was supported by the
German Research Foundation within the collaborative research center
TRR154 ``Mathematical Modeling, Simulation and Optimisation Using
the Example of Gas Networks'' (DFG-SFB TRR154/1-2014, TP B01), 
the Graduate School of Excellence Computational Engineering (DFG GSC233), 
and the Graduate School of Excellence Energy Science and 
Engineering (DFG GSC1070).

\bibliographystyle{plain}
\bibliography{bibimexpeer}

\end{document}